\documentclass[10pt]{article}

\usepackage{grffile}

\usepackage{amsmath,amsfonts,amscd,a4wide}
\usepackage{amsmath}
\usepackage{amssymb} 
\usepackage{graphicx}
\usepackage{amsthm}
\usepackage{color}
\usepackage{setspace} 

\usepackage{transparent}
\usepackage[small]{caption}
\usepackage{subcaption}
\usepackage{enumitem}

\usepackage{color,soul}
\definecolor{lightgray}{gray}{0.85}
\sethlcolor{lightgray}
\usepackage[margin=1in]{geometry}

\newtheorem{Lemma}{Lemma}[section]
\newtheorem{Theorem}{Theorem}
\newtheorem{Proposition}[Lemma]{Proposition}

\newtheorem{Remark}[Lemma]{Remark}

\newtheorem{Hypothesis}[Lemma]{Hypothesis}

\newenvironment{Proof}%
 {\begin{trivlist} \item[]{\bf Proof. }}%
 {\hspace*{\fill}$\rule{.4\baselineskip}{.4\baselineskip}$\end{trivlist}}

\newenvironment{Acknowledgment}%
 {\begin{trivlist}\item[]\textbf{Acknowledgments.}}{\end{trivlist}}

\makeatletter\@addtoreset{figure}{section}\makeatother

\makeatletter \@addtoreset{equation}{section} \makeatother

\newcommand{\R}{\mathbb{R}}
\newcommand{\C}{\mathbb{C}}

\newcommand{\Z}{\mathbb{Z}}
\newcommand{\mL}{\mathcal{L}}
\newcommand{\mF}{\mathcal{F}}
\newcommand{\mX}{\mathcal{X}}
\newcommand{\mY}{\mathcal{Y}}

	\newcommand{\mc}[1]{\mathcal{#1}}
	\newcommand{\mb}[1]{\mathbb{#1}}
	\newcommand{\tl}[1]{\tilde{#1}}
	\newcommand{\lp}{\left}
	\newcommand{\rp}{\right}
	\newcommand{\beq}{\begin{equation}}
	\newcommand{\eeq}{\end{equation}}
	\newcommand{\ba}{\begin{align}}
	\newcommand{\ea}{\end{align}}
	\newcommand{\fr}[2]{\frac{#1}{#2}}
	\newcommand{\p}{\partial}
	\newcommand{\ri}{\mathrm{i}}
	\newcommand{\rmi}{\mathrm{i}}
	\newcommand{\rlin}{\mathrm{lin}}
	\newcommand{\rmd}{\mathrm{d}}
	\newcommand{\re}{\mathrm{e}}

	\newcommand{\rmO}{\mathcal{O}}

	\newcommand{\rre}{\mathrm{Re}}
	\newcommand{\rim}{\mathrm{Im}}
	\newcommand{\rab}{\mathrm{abs}}

			\newcommand{\rs}{\mathrm{s}}
				
					\newcommand{\rcu}{\mathrm{cu}}
					\newcommand{\rcs}{\mathrm{cs}}
						\newcommand{\ru}{\mathrm{u}}
						\newcommand{\rc}{\mathrm{c}}

		\newcommand{\la}{\langle}
			\newcommand{\ra}{\rangle}

\newcommand{\br}{	\mathrm{br}}

\newcommand{\Addresses}{{
  \bigskip
  \footnotesize

R. ~Goh, \textsc{School of Mathematics, University of Minnesota,
    206 Church St. SE, Minneapolis, 55455 }\par\nopagebreak
  \textit{E-mail address}: \texttt{gohxx037@umn.edu},
  \textit{Phone}: \texttt{(612)-624-3531}

  \medskip

  A. Scheel, \textsc{School of Mathematics, University of Minnesota,
    206 Church St. SE, Minneapolis, 55455 }\par\nopagebreak
  \textit{E-mail address}: \texttt{scheel@umn.edu}
  \medskip

}}

\title{Hopf bifurcation from fronts in the Cahn-Hilliard equation}
\author{Ryan Goh and Arnd Scheel}

\begin{document}

\maketitle
\begin{abstract}

We study Hopf bifurcation from traveling-front solutions in the Cahn-Hilliard equation. The primary front is induced by a moving source term. Models of this form have been used to study a variety of physical phenomena, including pattern formation in chemical deposition and precipitation processes. Technically, we study bifurcation in the presence of essential spectrum. We contribute a simple and direct functional analytic method and determine bifurcation coefficients explicitly. Our approach uses exponential weights to recover Fredholm properties and spectral flow ideas to compute Fredholm indices. Simple mass conservation helps compensate for negative indices.  We also construct an explicit, prototypical example, prove the existence of a bifurcating front, and determine the direction of bifurcation.

\end{abstract}
\Addresses

\section{Introduction and main results}
\subsection{Motivation}

The Cahn-Hilliard equation 
$$
u_t = - ( u_{xx} + f(u))_{xx}, \quad f(u) = u - u^3
$$
was first proposed in \cite{Cahn58} to model the phase separation of a metal alloy under rapid homogeneous quenching.  Since then, it has been used to model a multitude of phase separation processes throughout the sciences; see \cite{fife02} or \cite{novick08} for an introduction and review of the equation and its properties. 
In particular, the Cahn-Hilliard equation has often been used to study pattern formation via chemical deposition and precipitation. Experiments studying such mechanisms date back to the time of Liesegang \cite{liese}, who studied the formation of periodic rings, now named after him, precipitating in the wake of a circular reaction front traveling through a gel solution.   These processes have since been found to create an incredible array of spatial patterns, ranging from regular structures such as periodic stripes and dot arrays, to more complex ones such as helices, chevrons, and fractals (see \cite{Droz00}, \cite{keller81}, \cite{Rajaram09},\cite{Lagzi13} and references therein). 

In the present, the quest continues to not only understand how such patterns arise, but also harness their power to create functional structures at the micro- and nanoscale; see for example  \cite{mijatovic05}, \cite{Squires05}, and \cite{voros05}.  In order to achieve regular patterns, one must control how and when instabilities are allowed to nucleate.  This can be achieved by using a triggering mechanism which travels through the medium, locally exciting the system as it travels. Along these lines, the Cahn-Hilliard model has been used to generally study such spatially progressive pattern formation via directional quenching fronts in  \cite{Foard12}, \cite{krekhov}.  More specifically, in chemical deposition and precipitation such triggering mechanisms typically take the form of a moving source which deposits mass, moving a stable medium into an unstable state. As the speed of this source varies, different patterns may be left in the wake. 

Specific examples of this type of deposition process arise in controlled evaporation, or "de-wetting" processes.  Here, a material is deposited onto a substrate through the spatially progressive evaporation of a solvent;  see \cite{Thiele14} for an in depth review of the many phenomena which can occur.  The Cahn-Hilliard equation has been used to model these phenomena with the variable $u$ representing the concentration of material deposited on the substrate.  In the work of \cite{Kopf12}, numerical continuation has been used to find modulated and un-modulated traveling wave solutions, revealing a rich snaking structure of saddle-node and Hopf-bifurcations as the speed of the evaporation front is varied.

\subsection{Our Setting}
Motivated by the aforementioned studies, we analyze equations of the form
\begin{equation}\label{e:CH}
u_t = - \lp( u_{yy} + f(y - ct,u)\rp)_{yy} + c\chi(y - ct;c)
\end{equation} 
Here, $u(y,t)$ is an order parameter which denotes the concentration of precipitate at a certain space-time point $(y,t)$ and $\chi(y- ct)$ is a source term which travels through the domain with a constant speed $c$, leaving behind a monotone, uniformly translating front $u_*(y,t) = u_*(y-ct)$.  The spatially dependent nonlinearity $f$ encodes any possible changes in the medium. 

For example, precipitation models such as \cite{Droz00} and \cite{Lagzi13} let $f(y-ct,u)\equiv f(u)$ and $\chi =  \chi(y - ct) $ be a localized gaussian source term. The resulting front $u_*(y-ct)$ then satisfies $u_*(y - ct)\rightarrow u_\pm$ as $y\rightarrow \pm \infty$ for each fixed $t$ and some constants $u_\pm$ with $u_+ - u_- = \int_{-\infty}^\infty \chi(\xi;c) d\xi$. 

Alternatively, the deposition models of \cite{Kopf12} and the directional quenching models of \cite{krekhov} have no source term, $\chi\equiv0$, and a nonlinearity $f$ which, in the co-moving frame $x = y-ct$, asymptotically approaches functions $f_\pm(u)$ as $x\rightarrow \pm \infty$ for all $x\in \R$.  In this case, patterns bifurcate from a trivial front $u_*(x)\equiv 0$.

\begin{figure}[h]

\includegraphics{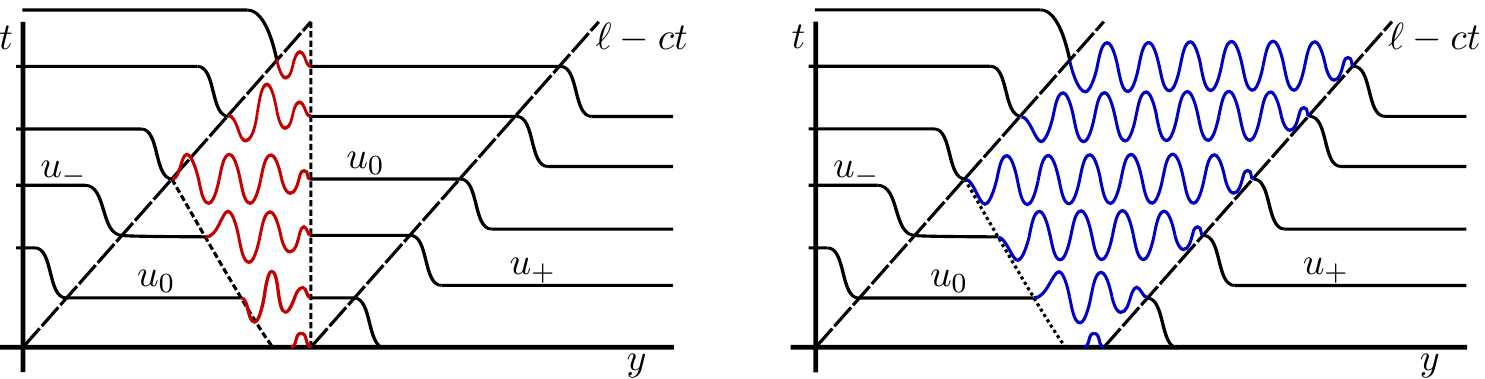}
  \caption{Development of an instability in spinodal regime of the traveling front $u_*$. Left: Instability is stationary (convective in the moving frame), and is "eaten" by the trailing homogeneous state. Right: Instability is absolute and is sustained as the front propagates through the medium.}
  \label{f:convabs}
\end{figure}

We study the behavior of the system near fronts $u_*(y-ct)$ which connect two homogenous equilibria lying outside or barely inside the spinodally unstable regime, where $f'>0$; see \cite{fife02}. As this front travels, there must be a moving spatial domain $[ -\ell - ct ,\ell - ct]$ where the front takes values inside this unstable regime. In the moving frame coordinate $x:=y - ct$, for large trigger speeds $c$, instabilities which may arise in this domain are convective (\cite{Scheel00}, \cite{Tobias98}), and get absorbed into the homogeneous equilibrium in the wake.  As $c$ decreases through a certain threshold, an absolute instability may arise, causing the formation of a periodic pattern which saturates the moving domain. In the physical literature, such a "self-sustaining" pattern is referred to as a nonlinear global mode; see \cite{Couairon97}.  See Figure \ref{f:convabs} for a schematic plot of these two types of instabilities. 

\subsection{Summary of contributions}
Our main contributions are threefold: a technically simple proof for Hopf bifurcation in the presence of essential spectrum, a rigorous existence proof of patterned fronts in the Cahn-Hilliard equation, and the analysis of an explicit, prototypical example. 
 
\paragraph{A simple proof of Hopf bifurcation in the presence of essential spectrum.}
When studying bifurcation in the presence of essential spectrum, the absence of a spectral gap precludes any center manifold reduction and hence does not allow the immediate use of standard finite-dimensional bifurcation techniques. Several methods have been developed to get around this difficulty.  In the context of Hopf bifurcation in viscous shocks and combustion waves, Texier and Zumbrun \cite{Zumbrun05,Texier08} use subtle point-wise estimates for the temporal semi-group of the linearization to obtain solutions via a Poincare return map. Alternatively, a more geometric approach via spatial dynamics can be used to study the spatial evolution of time-periodic functions. Bifurcating solutions then are constructed via a point-wise matching of infinite dimensional invariant manifolds. Such techniques have been applied in the study of viscous shocks by Sandstede and Scheel \cite{Sanstede08}, extending their previous work on Hopf bifurcations due to essential spectrum crossing the imaginary axis \cite{ssess1,
ssess2,ssess3}, and also in the propagation of water waves by Barrandon, Dias, and Iooss \cite{Barrandon05}, \cite{dias03}. Interaction between Hopf bifurcation and essential spectra was studied in a spirit similar to \cite{Zumbrun05,Texier08}, yet without the assumption of a conservation law, in \cite{kunze1,kunze2}, exploiting the sufficiently strong diffusive decay of modes associated with essential spectra in higher space dimensions. Hopf bifurcation in the presence of essential spectra is also responsible for the meandering transition of spiral waves \cite{barkley,ssw1,ssw2,gms}. While rigorous Hopf bifurcation results are not available for Archimedean spirals, the essential spectrum has striking consequences for the shape of bifurcating patterns, creating an intricate rotating super-spiral structure \cite{sssuper}.

Our main contribution is a much more direct functional analytic approach to this problem. While we develop our approach in the context of the Cahn-Hilliard equation, it could also be used in the context of \cite{Zumbrun05,sandstede08} to significantly simplify proofs. Our approach is simpler as it restricts the linear analysis to Fredholm properties on time-periodic functions, avoiding the subtle diffusive decay properties for infinite times used in \cite{Zumbrun05} or the point-wise reduction based on exponential dichotomies in \cite{sandstede08}.

In the setting described in the previous section, our approach exploits the  techniques in \cite{robbin95} to determine that the linearized equation is Fredholm with index -1 when considered on a suitable space of functions with spatial exponential weights and imposed  temporal periodicity. Mass conservation then allows us to restrict the codomain of the nonlinear operator to a certain subspace where its linearization has index 0. We then apply a Lyapunov-Schmidt reduction to this restriction to obtain existence of bifurcating solutions.

We also add that our method gives computable expressions for bifurcation coefficients.  In previous studies, such coefficients appear difficult to obtain; see  for example Eqn. 3.35 of \cite[\S 3.2]{Sanstede08}. Finally, we remark that this abstract approach should be applicable in many other situations, a few of which we discuss in Section \ref{s:dis} below.

\paragraph{Existence of pattern forming fronts.}
 
Our results show the existence of pattern forming fronts in the Cahn-Hilliard equation (\ref{e:CH}). As evidenced above, such fronts have been widely studied experimentally, numerically, and analytically. Furthermore, the computability of the bifurcation coefficients we obtain allows for the characterization of bifurcations and hopefully a deeper understanding of the patterns being formed.
 
\paragraph{Explicit characterization of a prototypical example.}
Finally, we apply our results to an idealization of the motivating examples discussed above which exhibits many interesting phenomena.  In particular, we study a nonlinearity of the form $f(y-ct,u) = \chi(y - ct) u +\gamma u^3  - \beta u^5$ and solutions which bifurcate from a trivial front $u_*\equiv 0$.  Here $\beta>0$, $\chi\equiv 1 $ for all  $x=y - ct \in[-l,l]$, and $\chi\equiv -1$ elsewhere. As it travels through the domain, the triggering mechanism $\chi$ does not add mass to the system but instead alters the stability of the homogeneous solution $u_*$. Indeed $\p_u f(x,0)>0$ (spinodally unstable) for all $x\in[-l,l]$ and $\p_uf(x,u_0)<0$ (spinodally stable) for all $x$ outside it.  As $c$ decreases through a certain speed $c_*$, we show that there exists a first-crossing of a pair of eigenvalues with non-zero imaginary part. The piecewise constant dependence of $f$ on $x$ allows us to determine leading order expansions for the accompanying eigenfunctions, for $l$ sufficiently large. We then apply 
our main result to conclude the existence of a bifurcating solution and furthermore that the bifurcation is subcritical for $\gamma>0$ and supercritical for $\gamma<0$.

\subsection{Hypothesis and main existence result}\label{ss:hmr}
In order to perform our analysis, we pass to a co-moving frame $x = y - ct$ so that \eqref{e:CH} becomes 
\begin{equation}\label{e:CH-mf}
u_t = -\lp[ u_{xx} + f(x,u) \rp]_{xx}   + c u_x+c\,\chi(x;c).
\end{equation}
We now specify the assumptions needed for our main result.

\paragraph{Nonlinearity and trigger.}
We start with assumptions on $f$ and $\chi$. 

\begin{Hypothesis}\label{h:nl}
The nonlinearity $f$ is smooth in both $x$ and $u$, and converges with an exponential rate to smooth functions $f_\pm := f_\pm(u)$ as $x\rightarrow \pm \infty$. This convergence is uniform for $u$ in bounded sets.
\end{Hypothesis}

\begin{Hypothesis}\label{h:us}
The trigger $\chi = \chi(x;c_*)$ is smooth and exponentially localized in $x$.
\end{Hypothesis}

\paragraph{Piecewise-smooth nonlinearity.}
Our explicit example, and several explicit models mentioned above are formulated in terms of discontinuous nonlinearities. We therefore include an alternate setup to cover such cases. 
\begin{Hypothesis}\label{h:up}
Let $l>0$, $\chi \equiv 0$, $u_*(x;c)\equiv u_-$ and $f(x,u) = b(x) (u - u_-) + g(u),$ where $b$ is piecewise smooth in $x$ with jump discontinuities at $x = \pm \ell$, and $g$ is smooth in $u$ such that $g(u) = \mc{O}(|u - u_-|^2)$. \end{Hypothesis}

\begin{Remark}\label{r:mjp}
Similar results will follow in the same manner if $b(x)$ has any finite number of jump discontinuities in $x$. For more general forms of $f$ which possess $x$-discontinuities that depend nonlinearly on $u$, our results should still hold but more complicated modifications to the smooth case are required.  
\end{Remark}

\paragraph{Existence and robustness of trigger front.}

We assume existence of a ``generic'' propagating front. 
\begin{Hypothesis}\label{h:extr}
There exists a front solution $u_*(y-c_*t;c_*)$ of \eqref{e:CH} for some $c_*>0$, with 
\[
\lim_{x\to\pm\infty}u_*(x;c_*)=u_\pm,\quad u_+-u_-=\int_\R u_*(\xi;c_*)\rmd\xi.
\]
Moreover, $u_*\in C^4(\R)$ and   
\[
|u_*(x) - u_\pm| +\sum_{j=1}^3|\p_x^{i} u_*(x)| \leq C'\re^{-\gamma|x|},
\]
for some $C,\gamma>0$. We refer to this front solution as the primary \emph{trigger front}
\end{Hypothesis}
One can show, under appropriate  assumptions on the nonlinearity $f$ and $u_\pm$, that such trigger fronts necessarily exist. One can indeed find those as solutions to a non-autonomous, three-dimensional ODE with a gradient-like structure, using Conley's index; see \cite{gohsch}.

We are interested in Hopf bifurcations from $u_*$. In the following state our spectral assumptions and their immediate consequences.

\begin{Hypothesis}\label{h:l0}
The point $0\in \C$ is not contained in the extended point spectrum of the linearization $L: H^4(\R)\subset L^2(\R)\rightarrow L^2(\R)$ defined as
\beq\label{e:Lx}
Lv:= -\p_{x}^2\lp(\p_{x}^2 v+ \p_uf(x,u_*(x)) v\rp) + c_* \p_x v.
\eeq
\end{Hypothesis}
\begin{Remark}\label{r:kerL}
This hypothesis implies that $\ker L = \varnothing$ when considered on the spaces $L^2(\R)$, $L_{\eta}^2(\R):= \{ u: \re^{\eta\sqrt{1+x^2}} u(x)\in L^2(\R)\}$, and $L_{+\eta}^2(\R):=\{ u: \re^{\eta x} u(x)\in L^2(\R)\}$ for any $\eta>0$ small.  It also follows that, when considered on the last of the spaces just listed, $L$ is invertible. For background on the notion of extended point spectrum, see for example \cite{fiesch}.
\end{Remark}
The following lemma guarantees that our assumptions so far are open in the class of problems considered here. 
\begin{Lemma}[Robustness of Front Solution]\label{l:robus}
Assuming the above hypotheses, for $\tl u_\pm$ in a small neighborhood of $u_\pm$, with $\tl u_+-\tl u_- = u_+-u_-$, and $c$ close to $c_*$, there exists a family of smooth front solutions $u_*(x;c)$ asymptotic to $\tl u_\pm$, satisfying Hypothesis \ref{h:extr}.
\end{Lemma}
\begin{Proof}
Hypothesis \ref{h:l0} implies that the steady-state equation associated with \eqref{e:CH-mf} (known as the traveling-wave equation) has a transverse intersection of the respective stable and unstable manifolds emanating from the hyperbolic equilibria $u\equiv u_\pm$. Indeed if this intersection was not transverse, the intersection would give rise to an exponentially localized solution of the linearized equation, hence contributing to the extended point spectrum.
\end{Proof}

\paragraph{Hopf crossing and non-resonance.}

We formulate our main spectral hypotheses on Hopf bifurcation.

\begin{Hypothesis}(simple Hopf-Crossing)\label{h:HC} 
The operator $L$, defined on $L^2(\R)$ as in \eqref{e:Lx} above, has a simple pair of algebraically simple eigenvalues $\lambda(c) = \mu(c) \pm\ri \kappa(c)$ and corresponding $L^2(\R)$-eigenfunctions $p(x), \overline{p}(x)$ such that for some $\omega_*\neq0$, and $c_*>0$ as above
$$
\mu(c_*) = 0,\quad \mu'(c_*) > 0,\quad\text{and}\quad \kappa(c_*) = \omega_*.
$$ 
\end{Hypothesis}
Note that the hypothesis implicitly assumes that $\rmi\omega_*$ does not belong to the essential spectrum, that is, $L-\rmi\omega$ is Fredholm with index 0. Let $\psi$ be the corresponding adjoint $L^2(\R)$-eigenfunction which is normalized so that
$$
\la \psi,p\ra_{L^2(\R)} = 1. 
$$  
Also, it can be readily obtained that 
$$
 \la \psi,\overline{p}\ra_{L^2(\R)} = 0.
 $$
Finally, we assume that there are no point or essential resonances:
\begin{Hypothesis}(Absence of resonances)\label{h:ex}
For all $\lambda\in \ri \omega_*\Z \,\diagdown \{0,\pm1\}$, the operator $L - \lambda$ is invertible when considered on the unweighted space $L^2(\R)$.
\end{Hypothesis}

\begin{Remark}\label{r:ex}
The Fredholm boundaries of $L - \lambda$ on the unweighted space $L^2(\R)$ are equal to the continuous curves
$$
\sigma_{\pm}:=\{\lambda\in\C\,: \, \lambda = k^4 - f_\pm'(u_\pm) k^2 - \ri c k,\, k\in \R\},
$$
Each of the curves $\sigma_\pm$ intersect the imaginary axis at $\lambda = 0$ and possibly two other points $\pm\ri \omega_e$. These last two intersections exist when $f'(u_\pm)>0$, respectively.
When considered on the doubly weighted space $L_\eta^2(\R)$ mentioned above,  the curves $\sigma_{\pm}$ are shifted
$$
\sigma^\eta_\pm:=\{\lambda\in\C\,: \, \lambda = (\ri k\mp\eta)^4 +  f_\pm'(u_\pm) (-k \mp \eta)^2 - c (\ri k \mp \eta),\, k\in \R\}.
$$
\end{Remark}

\begin{figure}[h!]
\includegraphics{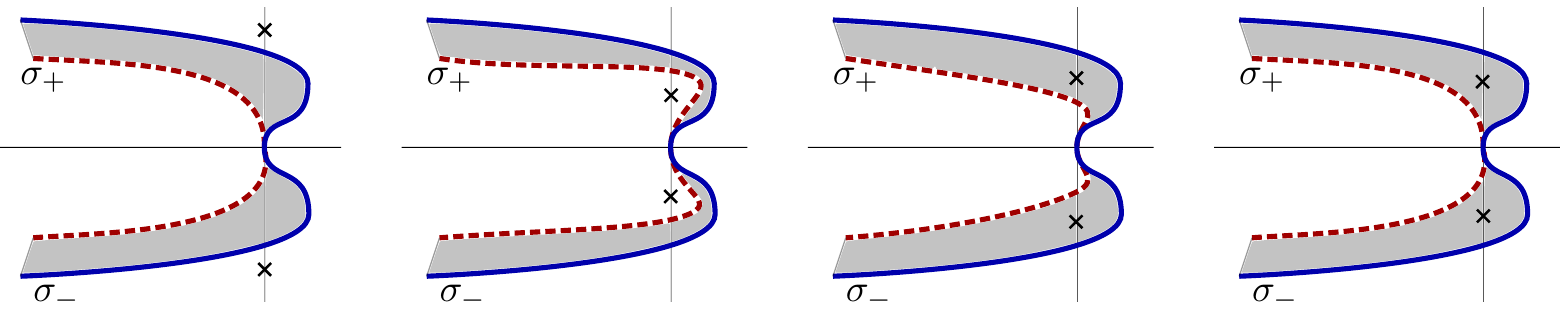}
  
  \caption{Examples of allowed (first two figures) and disallowed (last two figures) spectrum of $\mL$ in $\C$ under our hypothesis with $\eta = 0$.  The crosses denote the eigenvalues $\lambda_*(c_*), \overline{\lambda_*(c_*)}$, solid (blue) and dotted (red) lines denote the Fredholm borders $\sigma_\pm$ while the shaded region denotes the essential spectrum of $L$. (Color figure online)}
  \label{f:OK}
\end{figure}

Using the information in Remark \ref{r:ex}, Figure \ref{f:OK} depicts examples of Fredholm boundaries, for $\eta = 0$, which do and do not satisfy Hypothesis \ref{h:ex}.  The second figure from the left portrays the intriguing case where both $f'_\pm(u_\pm)>0$ and the Hopf eigenvalues are on the "wrong" side of the Fredholm borders $\sigma_\pm$. In other words, they lie inside of the essential spectrum of both of the constant coefficient operators $L_\pm:= -\p_{xx}\lp( \p_{xx} + f_\pm'(u_\pm)  \rp) + c_* \p_x$, but, since the Fredholm index is determined by the difference in Morse indices between $L_\pm$, $L - \lambda_*$ has index 0 and our spectral hypothesis are still satisfied.  Though our results give existence of time periodic solutions in this case, we believe that such solutions are not physically relevant: since the exponential weight selects the wrong spatial decay rates, the Hopf eigenvalues do not correspond to poles of the point wise Green's function.  Thus, any compactly supported initial data 
will be at most convectively unstable, leading to point wise decay as $t\rightarrow \infty$. Hence if the oscillatory part of our bifurcating  solution were multiplied by a compactly supported bump-function, it would decay as well for the linearized equation.

\paragraph{Large domain length.}  We may characterize the spectrum of $L$ more explicitly if we further restrict our hypotheses by assuming that $f(x,u)$ is piecewise constant in $x$ with $\p_u f(x,u_0)  \equiv C >0$ for all $x\in [-\ell,\ell]$. The results of \cite{gluing} imply that for $\ell>>1$ all but a finite set of the point spectrum of $L$ is well approximated by the absolute spectrum, $\Sigma_\rab\subset\C$, of the linearization about the homogeneous state $u_0$. Hence, eigenvalue crossings as described in Hypothesis \ref{h:HC} occur when $\Sigma_\rab$, typically through one or more branch points, crosses into the right half of the complex plane.  The front speed for which these crossings occur has come to be known as the linear spreading speed, which we denote as $c_\rlin$; see \cite{holzer14} and \cite{vanSaarloos03} for a more in depth discussion of these topics.  Thus, as $\ell\rightarrow \infty$, the Hopf-crossing speed $c_*$ will approach $c_\rlin$.  

In Section 4, our assumptions on $f$ allow us to use such an argument to prove the existence of a Hopf eigenvalue crossing and subsequently obtain explicit expansions for the corresponding eigenfunctions. This then allows us to apply Theorem \ref{t:hbex} to prove the existence of a Hopf bifurcation for an explicit example of this form.   

\subsection{Main Result}\label{s:mr}
We are now ready to state our main result.
\begin{Theorem}\label{t:hbex}
Given Hypotheses \ref{h:extr}, \ref{h:l0}, \ref{h:HC}, \ref{h:ex} and either the pair \ref{h:nl} and \ref{h:us}, or \ref{h:up}, there exists a one-parameter family of time-periodic solutions of \eqref{e:CH-mf} which bifurcate from the front solutions $u_*(x,c)$ as the speed $c$ decreases through $c_*$. This solution branch $(u,c,\omega)\in (u_*+H^4(\R))\times \R^2$  can be parameterized by $r\geq 0$, the amplitude of oscillations. More precisely, there exists $r_*>0$ and smooth functions $\Upsilon_j$, $j\in\{c,\omega,u\}$, defined for $|r|<r_*$, $\Upsilon_j(0)=0$, so that
\[
c=c_*+\Upsilon_c(r^2),\quad \omega=\omega_*+\Upsilon_\omega(r^2), \quad u=u_*+\Upsilon_u(r),
\]
with expansions
$$
\Upsilon_c(r^2)= \fr{\rre\{\theta_+\} }{\mu'(c_*)} r^2 + \mc{O}(r^4),\qquad
\Upsilon_\omega(r^2) = \rim\{ \theta_+\} |r|^2 + \mc{O}(|r|^4),\qquad \Upsilon_u(r)=r p \cos(\omega t)+\rmO(r^2).
$$
Here,  $\mu'(c_*)$ is the crossing speed from Hypothesis \ref{h:HC}, and 
$$
\theta_+=
\Big\la \lp(3\,\p_u^3f(x,u_*) \,p^2 \overline{p} \;+\;  \p_u^2f(x,u_*) \lp[p \varphi_0 + \overline{p} \varphi_{+}\rp]\rp)_{xx}, \psi\Big\ra_{L^2_\eta(\R)},
$$
with $p,\psi$ eigenfunctions and adjoint eigenfunctions, and   $\varphi_i$, defined in \eqref{e:LS3} below, encode quadratic interactions. In particular, if $\,\rre\{\theta_+\} >0$, the bifurcation is supercritical; if $\,\rre\{\theta_+\}<0$ then the bifurcation is subcritical. 
\end{Theorem}

\begin{figure}[h!]
        \centering
        \begin{subfigure}[b]{0.33\textwidth}
                \includegraphics[width=\textwidth]{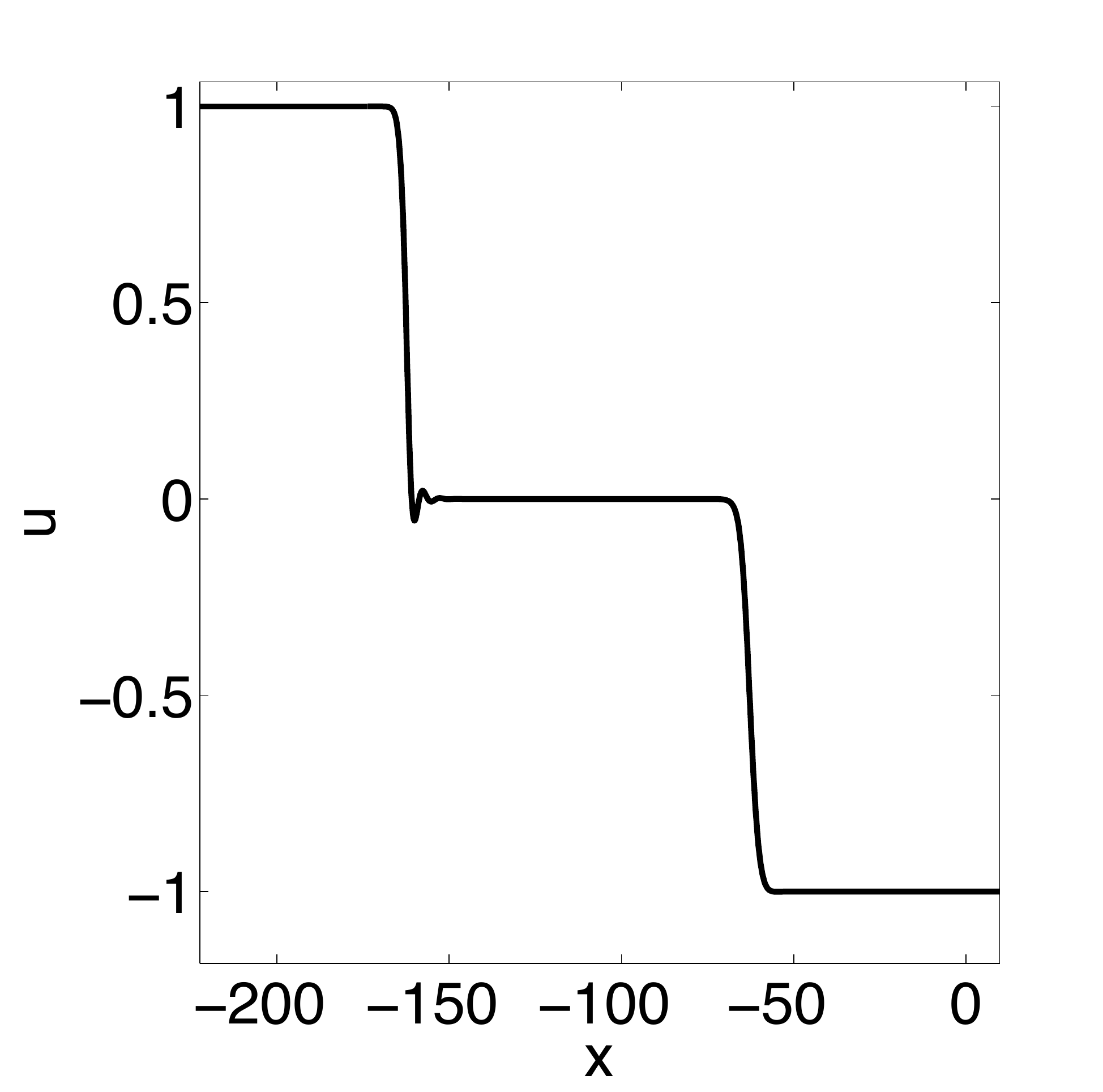}
                \caption{}
                \label{f:3a}
        \end{subfigure}%
        \begin{subfigure}[b]{0.33\textwidth}
                \includegraphics[width=\textwidth]{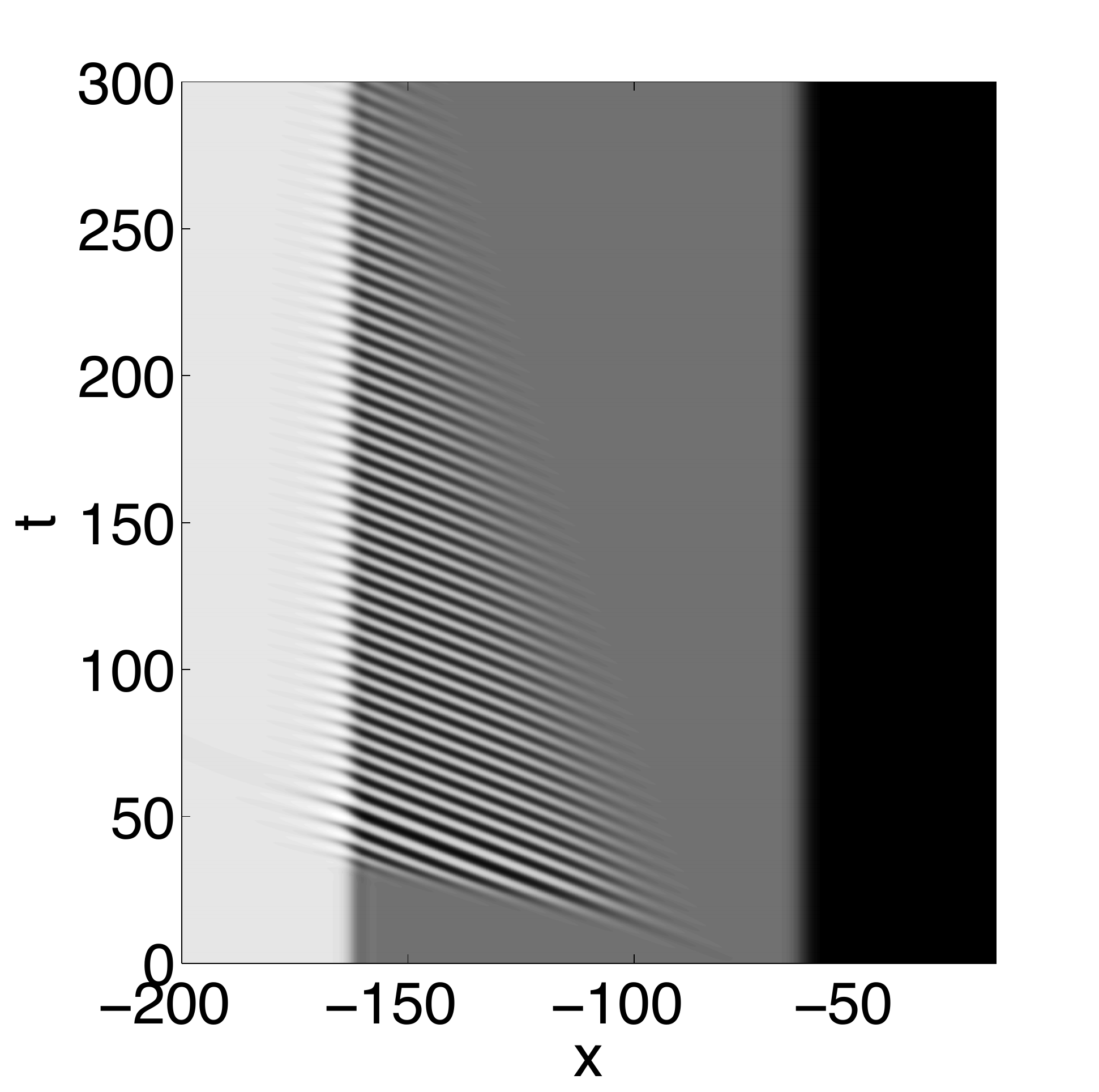}
                \caption{}
                \label{f:3b}
        \end{subfigure}
        \begin{subfigure}[b]{0.33\textwidth}
                \includegraphics[width=\textwidth]{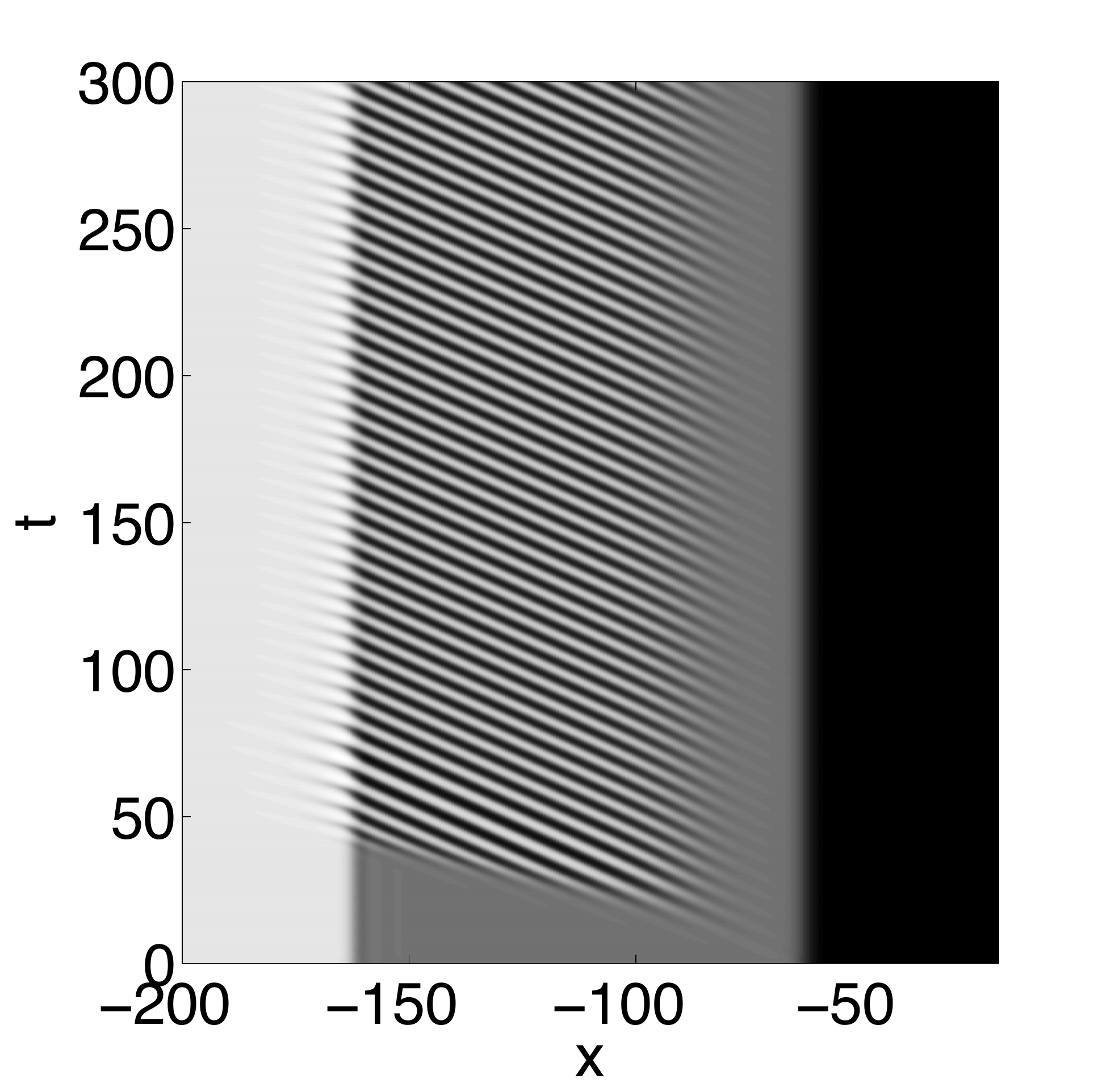}
                \caption{}
                \label{f:3c}
        \end{subfigure}
        \caption{(a): Front profile $u_*$ in co-moving frame for two Gaussian source terms. (b), (c): Spacetime plots in co-moving frame with speed $c>c_*$ and $c<c_*$ respectively. The initial condition for both is $u_*$ plus a small localized perturbation near $x = -75$.}\label{f:space-time}
\end{figure}

\begin{figure}[h!]
        \centering
        \begin{subfigure}[b]{0.33\textwidth}
                \includegraphics[width=\textwidth]{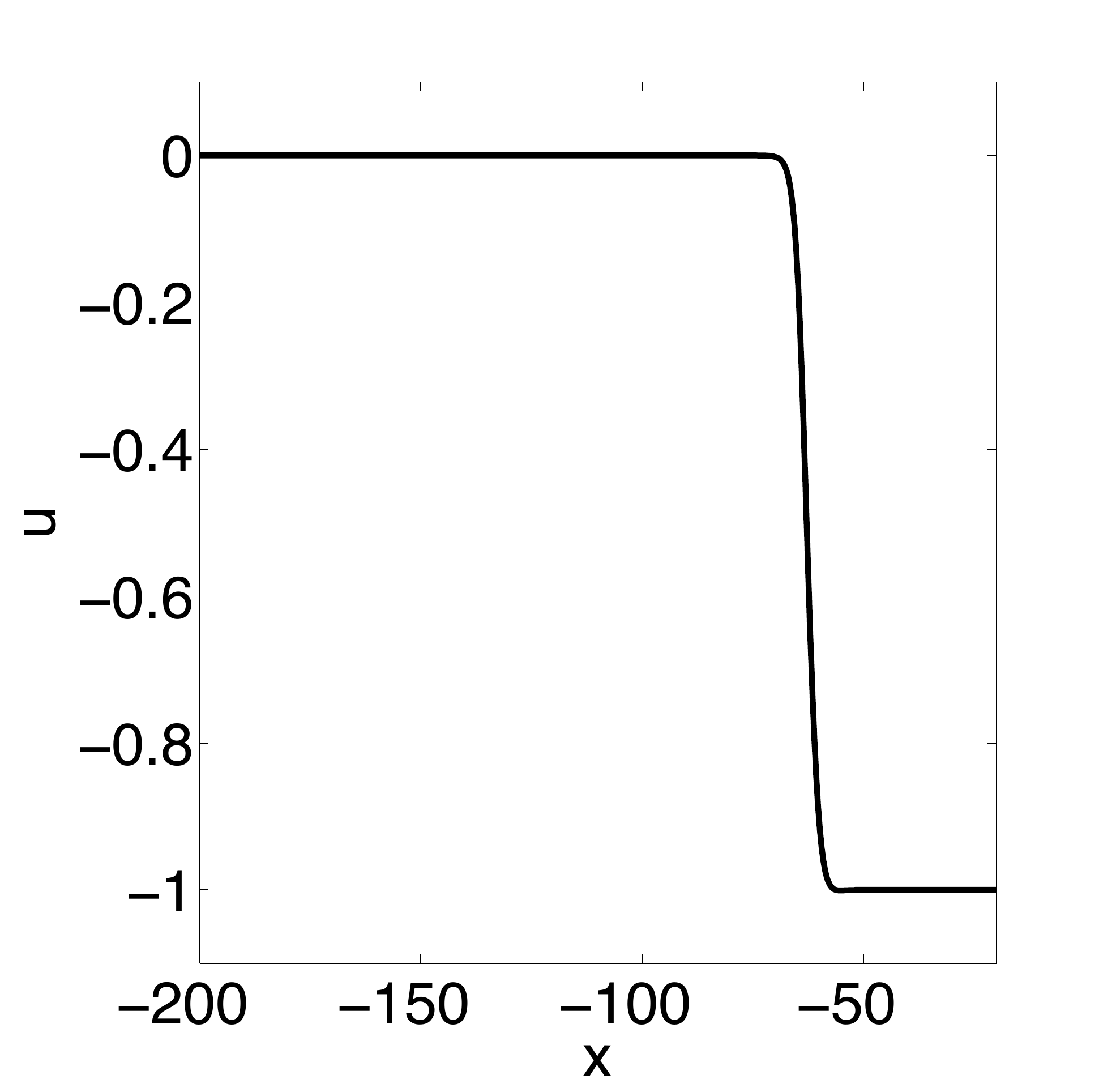}
                \caption{}
                \label{f:3a}
        \end{subfigure}%
        \begin{subfigure}[b]{0.33\textwidth}
                \includegraphics[width=\textwidth]{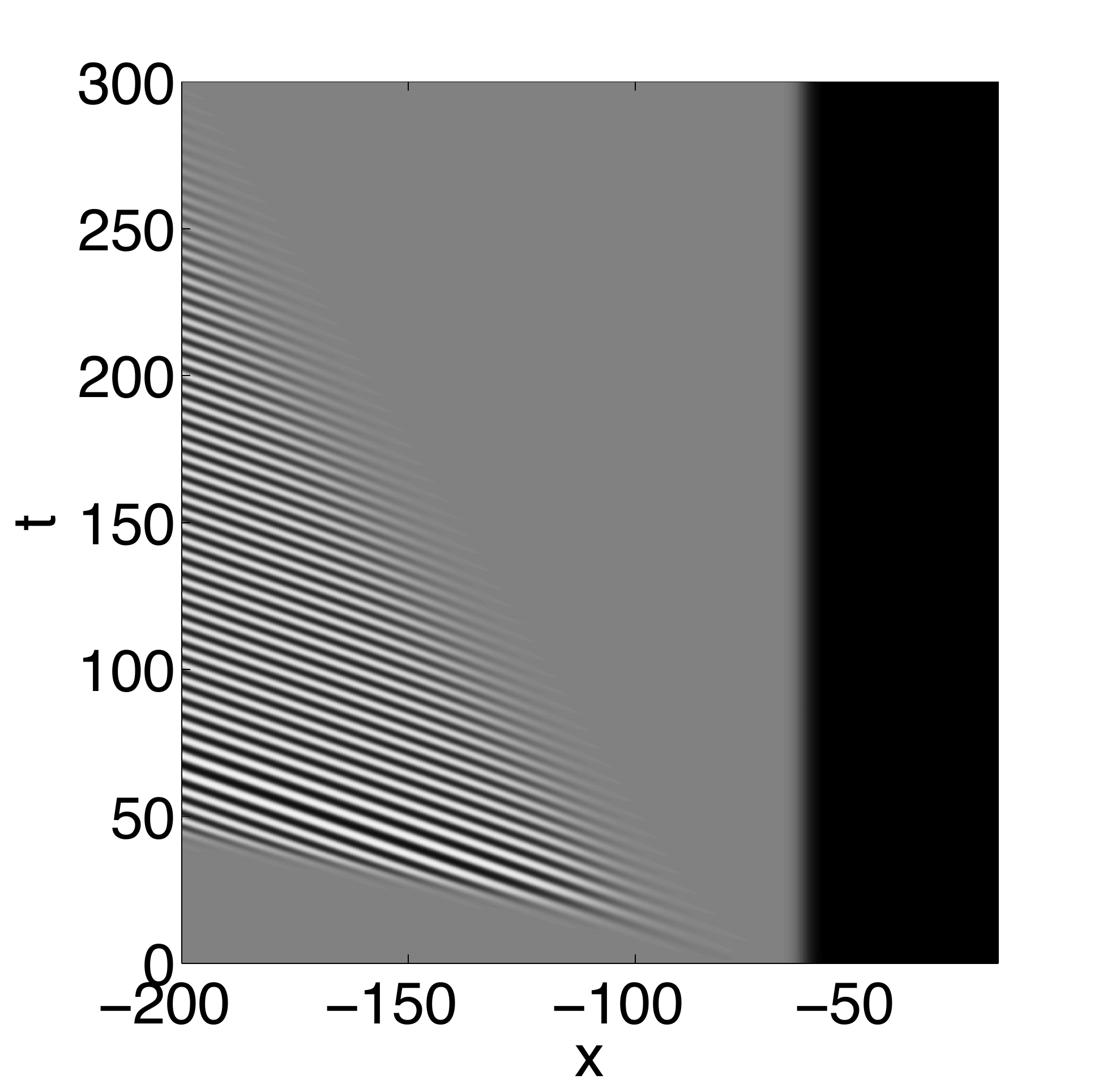}
                \caption{}
                \label{f:3b}
        \end{subfigure}
        \begin{subfigure}[b]{0.33\textwidth}
                \includegraphics[width=\textwidth]{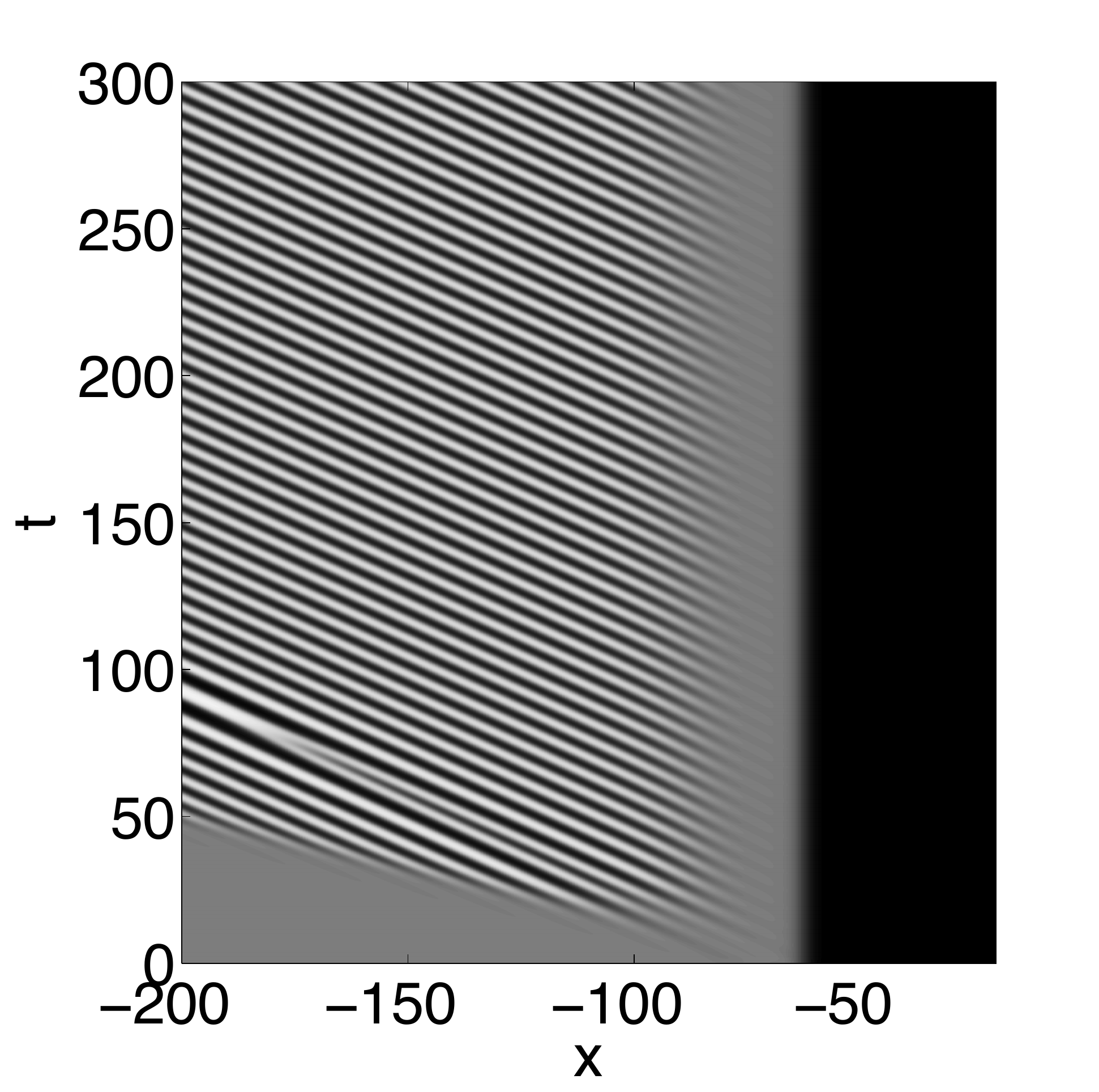}
                \caption{}
                \label{f:3c}
        \end{subfigure}
        \caption{(a): Front profile $u_*$ in co-moving frame. (b), (c): Spacetime plots in co-moving frame for speeds above and below a bifurcation point. The initial condition for both is $u_*$ plus a small Gaussian perturbation near $x = -75$.}\label{f:space-timeZ}
\end{figure}

A numerical example of this bifurcation is given in Figure \ref{f:space-time} where equation \eqref{e:CH-mf} is simulated with $f(x,u) = u - u^3$ and $\chi$ equal to a sum of two Gaussian source terms. The corresponding trigger front $u_*$ connects two stable homogeneous equilibria at $x = \mp\infty$ with a spinodally unstable plateau state in-between.  For speed $c>c_*$, oscillatory instabilities of this front are convected away, while for $c<c_*$ they are self-sustaining.  This setting, for which our results give a rigorous characterization, is closely related to models used by \cite{Droz00} and \cite{Lagzi13}, where $\chi$ is composed of only a single Gaussian and produces a front $u_*$ which connects a stable equilibrium at $x = +\infty$ to a spinodally unstable equilibrium at $x = -\infty$.  Here, a similar bifurcation occurs as the trigger speed is reduced.  Numerical simulations of such a situation are depicted in Figure \ref{f:space-timeZ}. 

\vspace{2in}

The remainder of the paper is organized as follows.  In Section \ref{s:pre}, we establish Fredholm properties of the linearization. This is done in Propositions \ref{p:fred} and \ref{p:fredJ}, for $f$ smooth and piecewise-smooth in $x$ respectively, using the methods of \cite{Faye13} and \cite{robbin95}. In Section \ref{s:pfT1} we give the proof of Theorem \ref{t:hbex}. In Section \ref{s:example} we study a prototypical example, showing the existence of a first crossing of Hopf eigenvalues. By finding leading order estimates for such a crossing and its corresponding eigenfunctions, we then apply Theorem \ref{t:hbex} to obtain existence of a Hopf bifurcation and compute the direction of branching. In Section \ref{s:dis} we discuss possible extensions and applications of our results.

\section{Preliminaries and Fredholm properties}\label{s:pre}

After introducing some notation and function spaces, we establish Fredholm properties in weighted spaces in Section \ref{ss:pfred}. We list necessary changes for the piecewise smooth case in Section \ref{ssec:disc}. 

For $\eta>0$ we define the exponentially weighted norm
\beq\label{e:wgt}
\| w\|_{ 2,\eta}^2 :=  \int_\R |\re^{\eta \la x\ra}  w(x)|^2 dx,
\eeq
where $\la x\ra  = \sqrt{1+x^2}$. We say that $w\in L^2_\eta(\mb{R})$ if $w$ is Lebesgue-measureable and $\|w\|_{2,\eta}<\infty$. Similarly, we define Sobolev spaces $H^k_\eta$, with $\partial^j u\in L^2_\eta$ for $j\leq k$. We use the following space-time norms, 
\begin{align}
X &= L^2(\mb{T}), \qquad Y = H^1(\mb{T}),\notag\\
\mc{X} &= L^2_{\eta}(\R, X),\qquad \mc{Y} = L^2_{\eta}(\R,Y) \cap H^4_{\eta}(\R,X),\label{d:sp}
\end{align}
where $x\in \R$ and $\tau \in \mb{T} = [0,2\pi)$, the one-dimensional torus.  Note that $Y$ is dense and compactly embedded in $X$. Also note that $\mX$ is a Hilbert space with the inner product
$$
\la u,v\ra_\mX: = \fr{1}{2\pi} \int_0^{2\pi} \int_{-\infty}^\infty u(x,\tau) \overline{v(x,\tau)} \,\,e^{2\eta\la x\ra}\,dx \,d\tau.
$$

Furthermore, the following norm makes $\mY$ a Hilbert space
\beq\label{e:nY}
\|u\|^2_\mY := \int_{-\infty}^\infty \|u(x,\cdot)\|^2_Y+ \sum_{j=1}^4\|\p_x^{j} u(x,\cdot)\|_X^2 dx.
\eeq

We define $\mc{F}:\mc{Y}\times \mb{R}^2 \rightarrow \mc{X}$ as
$$
\mF: (v;\omega, c)\longmapsto \omega v_\tau + \lp( v_{xx} + g(x,v; c)\rp)_{xx} - c \,v_x , \qquad  g(x,v; c) := f(x,u_*+v)  - f(x,u_*),
$$
so that time-periodic solutions $u = u_*+v$ of \eqref{e:CH-mf} satisfy $F(v;\omega,c) \equiv 0$.  

\subsection{Smooth nonlinearity}\label{ss:pfred}

We are interested in Fredholm properties of the linearization $\mL: \mY\subset \mX \rightarrow \mX$ of $\mF$ at the homogeneous solution $(v; \omega, c)  = (0; \omega_*,c_*)$, which has the form
\beq\label{e:L}
\mL : v\longmapsto \omega_* \p_\tau v - L v = \omega_*\p_\tau v - \p_x \tl L v,\qquad \tl L: u\longmapsto  -\p_x\lp( \p_x^2v + \p_uf(x,u_*) v\rp) + c_* v.
\eeq
These properties will be necessary to implement the Lyapunov-Schmidt reduction used in the proof of Theorem \ref{t:hbex}. We prove that $\mL$ is Fredholm in Proposition \ref{p:fred}  and we compute its Fredholm index in Proposition \ref{p:fInd}.

\begin{Proposition}\label{p:fred}
Assuming Hypotheses \ref{h:nl}, \ref{h:us}, \ref{h:l0}, \ref{h:HC}, and \ref{h:ex}, the operator $\mL:\mY\subset\mX\rightarrow \mX$ is Fredholm.
\end{Proposition}

Before proving the proposition, we prove the following lemma which adapts the methods of \cite{robbin95}; see also \cite{Faye13}. For $J>0$, let $\mY(J)$ and $\mX(J)$ denote the spaces of functions, in $\mY$ and $\mX$ respectively, which have $x$-support in the interval $[-J,J]$. Since the embedding $\mY(J)\hookrightarrow \mX(J)$ is compact, the following lemma allows us to apply an abstract closed range lemma \cite[Prop. 6.7]{taylor96}, showing that $\mL$ has closed range and finite dimensional kernel.

\begin{Lemma}\label{l:CR}
There exist constants $C>0$ and $J>0$ such that the operator $\mL$ defined above satisfies
\begin{equation}\label{e:CR}
 \|\xi\|_\mY  \leq C \lp(\|\xi\|_{\mX(J)} + \|\mL \xi\|_{\mX}\rp).
\end{equation}
\end{Lemma}
\begin{Proof}
Following \cite{robbin95}, the proof is divided into three steps:

\textbf{Step 1:}  Prove that the estimate holds for $J = \infty$.

For this step, momentarily assume that the exponential weight has $\eta = 0$.
To begin, we notice
\begin{equation}\label{e:CR1}
\|\mL \xi\|_{\mX} \geq \|(\p_\tau + \p_x^4)\xi\|_{\mX} - \|\p_x^2(\p_uf(x,u_*)\xi) - c_*\p_x\xi\|_{\mX}.
\end{equation}

Since $f$ and $u_*$ are smooth, for all $\epsilon>0$ we have
\begin{align}
\|\p_x^2(\p_uf(x,u_*)\xi) - c_*\p_x\xi\|_{\mX}&\leq C \|\xi\|_{H^2(\R,X)} \notag\\
&\leq C\|\xi\|_{\mX}^{\fr{1}{2}} \cdot \|\xi\|_{H^4(\R,X)}^{\fr{1}{2}}\notag\\
&\leq C\lp(\epsilon\|\xi\|_{H^4(\R,X)} + \fr{1}{4\epsilon} \|\xi\|_{\mX}\rp).\label{e:CR2}
\end{align}

Combining \eqref{e:CR1} and \eqref{e:CR2}, we have for sufficiently small $\epsilon>0$
\begin{align}\label{e:CR3}
\|\mL \xi\|_{\mX} + \fr{C}{4\epsilon} \|\xi \|_{\mX} &\geq \|(\p_\tau + \p_x^4)\xi\|_{\mX} - C\epsilon \|\xi\|_{H^4(\R,X)} \notag\\
&\geq C'\|u\|_{\mY}.
\end{align}
This gives the desired estimate.

If $\eta >0$ then this step works in essentially the same manner.  The only difference is that one must work with the conjugated operator $\mL_\eta:=\re^{\eta\langle x\rangle } \mL \,\re^{-\eta\langle x\rangle }$ and deal with third derivative terms which are small due to the fact that $\eta>0$ is small.

\textbf{Step 2}: Prove the estimate for the constant coefficient operators $\mL_\pm$ given above. 

We must work with the conjugated operators $\mL_{\pm,\eta}:=\re^{\eta\langle x\rangle } \mL_\pm \,\re^{-\eta\langle x\rangle } $.
By taking the Fourier transform in both $x$ and $\tau$, if 
\beq\label{e:cc}
\mL_{\pm,\eta} \xi = g,
\eeq
then
$$
\widehat{g}(\ri\zeta,\ri k) = \lp[(\zeta\mp \eta)^4 - f'_\pm(u_\pm)(\zeta\mp \eta)^2 - \ri c_*(\zeta\mp\eta) + \ri \omega k\rp]\widehat{\xi},\quad \xi\in \R, k\in\Z.
$$
By Hypothesis \ref{h:ex} and Remark \ref{r:ex}, for $\eta>0$ the essential spectrum of the time-independent operator does not intersect the set $\ri \omega_* \mb{Z}$.  Hence both equations in \eqref{e:cc} are invertible and 
$$
\widehat{\xi} =  \lp((\zeta\mp \eta)^4 - f'_\pm(u_\pm)(\zeta\mp \eta)^2 - \ri c_*(\zeta\mp\eta) + \ri \omega k\rp)^{-1} \widehat{g}
$$
such that the coefficient on the right is bounded.
The estimate
$$
\|\widehat{\xi}\|_\mX \leq \sup_{\zeta\in \R, k\in \Z} |\lp( (\zeta\mp \eta)^4 - f'_\pm(u_\pm)(\zeta\mp \eta)^2 - \ri c_*(\zeta\mp\eta) + \ri \omega k\rp)^{-1} | \cdot \|\widehat{g}\|_{\mX},
$$
implies by Fourier-Plancherel that 
$$
\|\xi\|_{\mY} \leq C_2 \|\mL_\pm \xi\|_{\mX}.
$$

\textbf{Step 3:}  Using estimates on $\mL - \mL_\pm$ which can be derived from Hypothesis \ref{h:us} and \ref{h:extr}, one performs a patching argument in the same way as in \cite{robbin95} (see also \cite{Faye13} for more details) to obtain the estimate \eqref{e:CR} for some sufficiently large $J>0$ and constant $C>0$.
\end{Proof}
\begin{Proof}[Proof of Prop. \ref{p:fred}]

Lemma \ref{l:CR} gives that $\mL$ has closed range and finite dimensional kernel.  To finish the proof we define a suitable adjoint $\mL^*$ and show it also satisfies a closed range lemma as above.  
In unweighted spaces the formal adjoint is
$$
\mL^*:= \p_x^4 +\p_uf(x,u_*) \p_x^2 +c_*\p_x -\omega\p_\tau.
$$
But as we wish to work with exponentially weighted spaces $\mY$ and $\mX$, we define the adjoint using the conjugated operator
$$
\mL^*_\eta := \re^{\eta\la x\ra} \mL^* \re^{-\eta\la x\ra}: H^4(\R,X)\cap L^2(\R,Y)\subset L^2(\R,X)\rightarrow  L^2(\R,X),
$$
Note, since $\mL_\eta$ is closed and densely defined, $\mL^*_\eta$ is well-defined.

This operator can then be run through the same estimates as in Lemma \ref{l:CR}, and we once again obtain that $\mL^*_\eta$ has closed range and finite kernel. Therefore $\mL$ is Fredholm. 

\end{Proof}
Next,  in Lemma \ref{l:i0} and \ref{l:ih}, we determine the index of $\mL$. We decompose $\mc{X}$ and $\mc{Y}$ into a direct sum of invariant subspaces so that $\mL$ is diagonal and the index of each restriction can be readily calculated.  Elementary results (see for example \cite{taylor96}) then give that $\mL$ has index equal to the sum of the indices of the restrictions.

The work of \cite[Thm. 1.5]{Arendt10} gives that $\mc{X}$ has a Fourier decomposition in the time variable
\begin{align}
\mc{X} &= \bigoplus_{k\in\mb{Z}} \mc{X}^k, \qquad \mc{X}^k := \lp\{ v\in \mc{X}: v(x,\tau) = \tl v(x) \re^{\ri k \tau}, \; \tl v\in L^2_\eta(\R) \rp\}.\notag
\end{align}

Next, let $\mc{X}_h = \bigoplus_{k\neq0} \mc{X}^k$, and  $\mc{Y}_h  = \bigoplus_{k \neq0} \mc{Y}^k,$ where $\mY^k= \mX^k\cap \mY$, so the following decompositions hold
$$
\mX = \mX^0\oplus \mX_h,\qquad \mY = \mY^0\oplus \mY_h.
$$ 
Note that $\mc{X}_0$ is the set of all time-independent functions in $\mc{X}$ while $\mX_h$ is the set of all functions with time-average equal to zero.  

\begin{Lemma}\label{l:i0}
The restriction $\mc{L}_0:= \mc{L}: \mc{Y}^0 \rightarrow \mc{X}^0$ has Fredholm index $-1$.
\end{Lemma}

\begin{Proof}
On $\mc{Y}_0$ we have $\mL = -L  = -\p_x\circ \tl L$.  Recall from \eqref{e:L} that $\tl L: H^3_\eta(\R)\subset L^2_\eta(\R)\rightarrow L^2_\eta(\R)$ is defined as $\tl L v =  -\p_x(\p_x^2 v+ \p_uf(x,u_*)v) + c_*v$.  Hypothesis \ref{h:us} implies that $\tl L$ is an asymptotically constant operator:
$$
\tl L \rightarrow \tl L_\pm \quad \text{as} \;\;x\rightarrow \pm\infty,\qquad \tl L_\pm v\,:=-\, \p_x^3 v -  f'_\pm(u_\pm)\p_xv +c_*v.
$$

Moreover, the constant coefficient first order systems associated with $\tl L_\pm$ are hyperbolic with the same Morse index. Indeed, since $c>0$, each of the polynomials $\nu^3+f'_\pm(u_\pm) \nu - c = 0$ has two positive roots and one negative root.  Thus, the piecewise constant operators $\tl L_+$ and $\tl L_-$ have  relative morse index equal to zero as well. This implies that the operator $\tl L$ has Fredholm index equal to zero. (See for example \cite[Sec. 3.1.10 - 11]{kapitula13})).   

To finish the proof it suffices to notice that $\p_x: H^1_\eta(\R)\rightarrow L^2_\eta(\R)$ is Fredholm with index -1.  The result then follows using standard results on the composition of Fredholm operators (see for example \cite[Sec. A.7]{taylor96})
\end{Proof}

\begin{Lemma}\label{l:ih}
The restriction $\mc{L}_h:=\mc{L}: \mc{Y}_h \rightarrow X_h $ has Fredholm index $0$.
\end{Lemma} 
\begin{Proof}
First note that since $\mL$ is Fredholm, $\mL_h$ and $\mL_h^*$ must have finite dimensional kernel.
Next, it is straightforward to notice that each restriction $\mL_k := \mL: \mc{Y}^k \subset \mc{X}^k \rightarrow \mc{X}^k$ is well defined and takes the form
$$
\mL_k (\re^{\ri k \tau}\tl v(x)) = \p_{x}^2\lp(\p_x^2\tl v + \p_uf(x,u_*(x))\tl v\rp) - c_* \p_x \tl v + \ri \omega k \tl v.
$$
Hypothesis \ref{h:l0}  and \ref{h:ex} imply that for $\mL_k$ and its adjoint $\mL_k^*$
$$\dim \ker \mL_k =
\begin{cases}
 &0\quad k\neq \pm 1\\
 &1 \quad k = \pm 1,
 \end{cases}  
 \qquad \dim \ker \mL_k ^*=
\begin{cases}
 &0\quad k\neq \pm 1\\
 &1 \quad k = \pm 1.
 \end{cases}
 $$
Comparing the dimensions of $\ker \mL_h$ and $\ker \mL_h^*$ then implies that $\mL_h$ has Fredholm index 0. 
\end{Proof}

\begin{Remark} For simplicity, we have included direct proofs to determine the Fredholm index of $\mL$ in the preceding lemmas. We note that one could also calculate these using a spectral flow as in \cite{robbin95}. Namely, the index could be found by tracking spatial eigenvalue crossings as $x$ moves from $-\infty$ to $+\infty$. 
\end{Remark}

The previous lemmas then give the following proposition.

\begin{Proposition}\label{p:fInd}
Given the hypotheses in Proposition \ref{p:fred}, the operator $\mL$ has Fredholm index -1.
\end{Proposition}
\begin{Proof}
Since $\mL$ can be decomposed as
$$
\mL:=\left(\begin{array}{cc}\mL_0 & 0 \\0 & \mL_h\end{array}\right): \mc{Y}_0\oplus \mc{Y}_h \longrightarrow \mc{X}_0\oplus \mc{X}_h,
$$
a standard result in Fredholm theory (see \cite[Sec. A.7]{taylor96}) gives that the Fredholm index of $\mL$ is equal to the sum of the indices of $\mL_0$ and $\mL_h$.  This fact, in combination with Lemmas \ref{l:i0} and \ref{l:ih}, proves the proposition.
\end{Proof}

\subsection{Piecewise smooth $f$}\label{ssec:disc}
If Hypothesis \ref{h:up} is assumed instead of Hypotheses \ref{h:nl} and \ref{h:us} the setting must be slightly altered in order to obtain the Fredholm properties required in the proof of Theorem \ref{t:hbex}.  In particular, since $\p_u f$ has discontinuities in $x$, $\mL$ is not well-defined on $\mY$ and hence jump-conditions are needed. Let us define the jump condition notation
$$
\delta_{x_0} u = \lim_{x\rightarrow x_0^+}u(x,t) -  \lim_{x\rightarrow x_0^-}u(x,t).
$$  
Also, for simplicity let us define the piecewise-smooth function 
$$
b(x) = \p_u f(x,u_*(x)).
$$
Next let us define the following set of conditions on a function $u(x,t)$,
\begin{equation}\label{e:jp}
(\#\#):= \begin{cases}
&t\in \mb{T},\,x_0 = \pm \ell,\notag\\
 &\delta_{x_0} u = 0,\quad \delta_{x_0} u_x = 0,\quad\delta_{x_0} u_{xx} = -u(x_0,t) \,\delta_{x_0}b ,\notag\\
 &\quad\delta_{x_0} u_{xxx}  = -\lp[u(x_0,t) \,\delta_{x_0}b_x + u_x(x_0,t)\,\delta_{x_0}b \rp].\notag
\end{cases}
\end{equation}

A brief calculation then shows that $\mL$ is well defined on the space
$$
\mY^{\#\#}:=  \Big( H^4_\eta\big([\ell,\infty), X\big) \oplus H^4_\eta\big([-\ell,\ell]), X\big)\oplus H^4_\eta\big((-\infty,-\ell], X\big) \Big) \cap L^2\big(\R,Y\big) \cap (\#\#) \subset \mX.
$$

Furthermore, it is easily seen that $\mL$ is a closed, densely defined operator on $\mX$. Indeed,  the latter fact follows from the density of $\mY$ in $\mX$, and the fact that for any $u\in\mY$ there exists a function $v$, which is smooth away from the points $x = \pm \ell$, has arbitrarily small $L^2$-norm, and yet satisfies the jump conditions $(\#\#)$ so that $u+v\in \mY^{\#\#}$. 
 
For ease of notation, we restrict for the remainder of the section to $f$ with one jump discontinuity located at $x = 0$. The result for a nonlinearity with multiple discontinuities will follow in the same manner. Hence we work with the operator
\beq
\mL:\mY^{\#}:=\lp(H^4_\eta(\R^-,X)\oplus H^4_\eta(\R^+,X) \rp) \cap L^2_\eta(\R,Y)\cap (\#)\subset \mX\rightarrow \mX,
\eeq
defined as in \eqref{e:Lx} above with 
\begin{equation}\label{e:jp0}
(\#):= \begin{cases}
\, &t\in \mb{T},\,x_0 =0\notag\\
 &\delta_{x_0} u = 0,\quad \delta_{x_0} u_x = 0,\quad\delta_{x_0} u_{xx} = -u(x_0,t) \,\delta_{x_0}b ,\notag\\
 &\quad\delta_{x_0} u_{xxx}  = -\lp[u(x_0,t) \,\delta_{x_0}b_x + u_x(x_0,t)\,\delta_{x_0}b \rp].\notag
\end{cases}
\end{equation}

Our approach is to conjugate from $\mY$ to $\mY^\#$ through a change of variables $\tl u = u + \Phi $ where $\Phi = \Phi(x,\tau)$ has jump discontinuities on $(0,\tau)\in \R\times \mb{T}$ which compensate for the discontinuities created by $b(x)$. We construct $\Phi$ using solutions of fractional order, $L^2(\mb{T})$-valued, evolution equations which have the jump conditions on $\{0\}\times \mb{T}$ as initial conditions.

For any $\mu,\tau\in [0,\infty)$, and open $U\subset \R$ let
$$
W^{\sigma,\mu}_\eta(U\times \mb{T}) := H^\mu_\eta(U, X)\cap L^2_\eta(\R, H^\sigma(\mb{T})),
$$
denote the anisotropic Sobolev space of order $\mu$ in space and order $\sigma$ in time defined in the usual way via Fourier Transform. As they will repeatedly arise in the following, denote $V = \R\times \mb{T}$ and $V^\pm = \R^\pm \times\mb{T}$. Note also that $W_\eta^{0,0}(V) = \mX$, $W_\eta^{1,0}(V) = L^2_\eta(\R,Y)$, $W_\eta^{0,4}(V) = H^4_\eta(\R,X)$, and $ W_\eta^{1,4}(V) = \mY$.

To setup the evolution equations, for $i = 0,1$ and $\alpha_i\in \R^+$, let $A_i: H^{\alpha_i}(\mb{T})\rightarrow L^2(\mb{T})$ be the linear operators defined via Fourier series as
\beq\label{e:fr-op}
(\widehat{A_i v})_k:= (-|k|^{\alpha_i}-1)\; \widehat{v}_k,\qquad k\in \Z.
\eeq
Next, for $\beta_i\in \R^+$, define the trace operators $T_i:W^{1,4}_\eta(V^+)\rightarrow H^{\beta_i}(\mb{T})$ as$$
T_0[u]:= -\delta_0(b)u(0,\tau),\qquad T_1[u]:=-(\delta_0(b_x) u(0,\tau) + \delta_0(b)u_x(0,\tau)).
$$ 
Anisotropic trace estimates give that $u(0,\tau)\in H^{7/8}(\mb{T})$ and $u_x(0,\tau)\in H^{5/8}(\mb{T})$ if $u\in W^{1,4}_\eta(V^+)$; see \cite[Lem 3.5]{denk05}. This means that $T_0$ and $T_1$ are well defined for $\beta_0\leq7/8$ and $\beta_1\leq 5/8$ respectively. Also, note that if these inequalities are strict then each $T_i$ is compact.

In order to obtain the desired regularity properties we set
\begin{equation}
\alpha_0 = \fr{5}{16}, \quad \alpha_1 = 1,\quad \beta_0 = \fr{7}{8} - \epsilon, \quad\beta_1 = \fr{9}{8} - \epsilon,
\end{equation}
for some $\epsilon>0$ sufficiently small. Then, define the $X$-valued initial value problems
\begin{align}
\p_xv_0 &= A_0 v_0,  \quad v_0(0) = T_0 [u]\label{e:fr-ev0},\\
\p_xv_1 &= A_0 v_1, \quad v_1(0) = T_1 [u] - A_0 T_0[u]\label{e:fr-ev1}.
\end{align}
where $v_i = v_i(x)$ take values in $X = L^2(\mb{T})$.  We then obtain the following result characterizing solutions of these equations.

\begin{Proposition}\label{p:ev-eq}
Given $u\in W^{1,4}_\eta(V^+)$, there exist unique solutions $v_{0}^*$ and $v_1^*$ of the initial value problems \eqref{e:fr-ev0} and \eqref{e:fr-ev1} which satisfy
\begin{align}
v_0^*&\in W^{\beta_0 - \alpha_0/2,1}_\eta(V^+) \cap W^{\beta_0 + \alpha_0/2,0}_\eta(V^+)\cap W^{\beta_0 - 3\alpha_0/2}_\eta(V^+),\\
v_1^*&\in W^{\beta_1 - \alpha_1/2,1}_\eta(V^+) \cap W^{\beta_1+\alpha_1/2,0}_\eta(V^+).
\end{align}
\end{Proposition}
\begin{Proof}
This result can be proved using Fourier analysis. For a more general reference see \cite{Amann95}
\end{Proof} 
Note that for the specific values of $\alpha_i,\beta_i$ listed in \eqref{e:fr-ev0} and \eqref{e:fr-ev1}, we have $v_0^*\in W^{1,2}_\eta(V^+)$ and $v_1^*\in W^{1,1}_\eta(V^+)$.  We may then define functions $\Phi_i = \Phi_i(x,\tau)$ as
\beq
\Phi_0(x,\tau) := \int_0^x\int^y_0 v_0^*(s,\tau)ds\,dy,\quad \Phi_1(x,\tau):=  \int_0^x\int^y_0\int_0^z v_1^*(s,\tau)ds\,dz\,dy,\quad
\eeq
so that $\Phi_0,\Phi_1\in W^{1,4}_\eta(V^+)$.  By extending $\Phi_i$ by zero for $(x,\tau)\in V^-$ and using the fact that 
$$(\Phi_0)_x\,\Big|_{x=0}  = (v_0^*)_x\,\Big|_{x=0} = A_0 v_0^*\,\Big|_{x=0}  = A_0T_0[u],$$
our construction gives that $\lp(\Phi_0 + \Phi_1\rp)$ satisfies the jump conditions $(\#)$ defined above. Hence the following mapping is well defined
\begin{align}
\Phi:W^{1,4}_\eta&(V)\rightarrow \mY^\#,\\
&u\longmapsto\Phi[u] = \Phi[u](x,\tau):= \rho(x)(\Phi_0(x,\tau)+\Phi_1(x,\tau)),
\end{align}
where $\rho = \rho(x)$ is a smooth bump function compactly supported and identically equal to 1 in a neighborhood of the origin. 
We then have the following lemma
\begin{Lemma}\label{l:iso}
The mapping $\mathrm{id} + \Phi: \mY\rightarrow \mY^\#$ is a linear isomorphism.
\end{Lemma}
\begin{Proof}
It is readily found that this mapping is linear.  Furthermore, since we have not used the full trace regularity of $u$, the mapping $\Phi$ is compact.  Hence, it suffices to show that $\mathrm{id} + \Phi$ is one-to-one, since it then immediately follows that the mapping is onto.  If $(u + \Phi[u]) = 0$, then for all $\tau\in \mb{T}$
$$
\tl u(0,\tau) = \Phi[u](0,\tau) = 0 = \Phi[u]_x(0,\tau) = u_x(0,\tau).
$$
This implies 
$$
\delta_0 (\Phi[u]_{xx} )= \delta_0(\Phi[u]_{xxx}) = 0,
$$
so that the initial value problems \eqref{e:fr-ev0} and \eqref{e:fr-ev1} have zero initial conditions and hence that $u = 0$.
\end{Proof}

We are now ready to prove the desired result. 
\begin{Proposition}\label{p:fredJ}
Assuming the Hypotheses \ref{h:l0}, \ref{h:HC}, \ref{h:ex}, and \ref{h:up},  the operator $\mL:\mY^\#\subset \mX \rightarrow \mX$ is Fredholm with index -1.
\end{Proposition}
\begin{Proof}

First note that because $\mL$ is closed and densely defined, its $\mX$-adjoint is
\begin{align}
\mL^*&: \mY \subset \mX \rightarrow \mX\notag\\
\mL^*v &:= -\p_t u -\p_x^4 u - b(x) \p_x^2u - c_* \p_x u.\label{e:LaP}
\end{align}
This definition can be easily calculated using the jump conditions in (\#\#) and integration by parts.  The methods used to prove Proposition \ref{p:fred} can immediately be applied to obtain that $\mL^*$ has closed range and finite dimensional kernel.

Since $\mathrm{id}+\Phi$ is an isomorphism, it suffices to prove that $\tl\mL := \mL\circ(\mathrm{id}+\Phi): \mY\rightarrow \mX$ has closed range and finite dimensional kernel. We thus proceed as in Lemma \ref{l:CR}.  We only give the proof of Step 1, obtaining a G\r{a}rding type inequality as in \eqref{e:CR3}. The subsequent steps will then follow in an analogous way to those in Lemma \ref{l:CR}. In particular, since $\tl\mL$ is equal to constant coefficient operators $\mL_\pm$ for $x$ outside the support of $\Phi[u]$, an even simpler patching argument than that of \cite{Faye13} and \cite{robbin95} can be implemented. Also, since we have conjugated to the space $\mY$, we still have the compact embedding of the truncated spaces $\mY(J)\hookrightarrow \mX(J)$ and may apply the abstract closed range lemma. Therefore $\tl\mL$ and, by Lemma \ref{l:iso}, $\mL$ have closed range and finite dimensional kernel.  Since $\mL^*$ has the same properties, we find that $\mL$ is Fredholm.  The index can be found in the exact same manner as in Section \ref{ss:pfred}.

Let $\lesssim$ and $\gtrsim$ denote inequality up to a constant independent of the variables being used. We first estimate
\begin{align}
\|\mL \Phi[u]\|_{\mX} &= \|\mL \Phi[u]\|_{W^{0,0}_\eta(V^+)}\notag\\
& \lesssim \|\Phi[u]\|_{W^{1,4}_\eta(V^+)} \lesssim \|\Phi[u]\|_{W^{1,0}_\eta(V^+)} + \|\Phi[u]\|_{W^{0,4}_\eta(V^+)}\notag\\
&\lesssim\|u\|_{W^{1-\epsilon,0}_\eta(V)}+ \|u\|_{W^{0,2}_\eta(V)}  \notag\\
&\leq  c_1(\epsilon)\|u\|_{\mX}+ c_2(\epsilon)\|u\|_{\mY}. \label{e:est1}
\end{align}
where $c_1(\epsilon)\rightarrow \infty$ and $c_2(\epsilon)\rightarrow 0$ as $\epsilon\rightarrow 0^+$. In the first two lines we restricted to $V^+$ because $\mathrm{supp}(\Phi)\subset V^+$. In the third line, the first term is obtained by exploiting the fact that less than maximal regularity of the trace, $u|_{x = 0}$, is used. The second term in the third line is obtained using trace and inverse-trace estimates from \cite[Lem. 3.5]{denk05}, 
\begin{align}\|\Phi[u]\|_{W^{0,4}_\eta(V^+)} &\lesssim \|\Phi[u]_{xx}\|_{W^{0,2}_\eta(V^+)}\lesssim \|u(0,\cdot)\|_{W^{3/4,3/2}_\eta(\{0\}\times \mb{T})} \notag\\
&\lesssim \|u\|_{W^{0,2}_\eta(V^+)} + \|u\|_{W^{0,2}_\eta(V^-)}\notag\\
&\sim \|u\|_{W^{0,2}(V)}.
\end{align}
To obtain the last line in \eqref{e:est1}, we use the estimates
\begin{align}
\|u\|_{W^{0,2}_\eta(V)} &\lesssim \fr{1}{\epsilon}\|u\|_{\mX} + \epsilon\|u\|_{\mY},\label{e:est2a}\\
\|u\|_{W^{1-\epsilon,0}_\eta(V^+)} &\leq \|u\|^\epsilon_{W^{0,0}_\eta(V)} \|u\|^{1-\epsilon}_{W^{1,0}_\eta(V)}\leq \tl c_1(\epsilon) \| u\|_{W^{0,0}_\eta(V))}+ \tl c_2(\epsilon)\|u\|_{W^{1,0}_\eta(V)},\label{e:est2b}
\end{align}
where $\tl c_1(\epsilon)\rightarrow+\infty$ and $\tl c_2(\epsilon)\rightarrow 0$ as $\epsilon\rightarrow 0$.
The estimate \eqref{e:est2b} uses standard Sobolev interpolation results (see \cite{lofstrom76}) and Young's inequality.

This finally allows us to obtain
\begin{align}
\| \tl\mL u\|_\mX& = \| \mL\circ(\mathrm{id} +\Phi) u\|_{L^2_\eta(V)}  \gtrsim \sum_{i =\pm}\| \mL\circ(\mathrm{id} + \Phi) u\|_{L^2_\eta(V^i)}\notag\\
 &\geq \sum_{i = \pm}\| \mL u\|_{L^2_\eta(V^i)} - \|\mL\Phi[u]\|_{L^2_\eta(V^i)} \notag\\
&\gtrsim \sum_{i = \pm}\lp(\|u\|_{W^{1,4}_\eta(V^i)} - \|u\|_{L^2_\eta(V^i)}\rp) - \| \tl\mL \Phi[u]\|_\mX\notag\\
&\gtrsim \sum_{i = \pm}\lp(\|u\|_{W^{1,4}_\eta(V^i)} - \|u\|_{L^2_\eta(V^i)}\rp) - \lp(c_1(\epsilon)\|u\|_{\mX}+ c_2(\epsilon)\|u\|_{\mY}\rp)\notag\\
&\geq C_\epsilon\|u\|_{\mY}  - C'_\epsilon\|u\|_{\mX},\label{e:est}
\end{align}
with $ C_\epsilon,  C'_\epsilon>0$ for $\epsilon>0$ sufficiently small. Since $\mL\circ(\mathrm{id} - \Phi) u\in L^2_\eta(V)$, the first and last inequality follow from the equivalence of the Euclidean and box norms on $\R^2$. The third inequality is obtained by proceeding as in Step 1 of Lemma \ref{l:CR} on each $V^i$.  The estimate \eqref{e:est1} gives the fourth inequality.
\end{Proof}
In the more general case where $b(x)$ has more than one discontinuity, one must first construct jump functions $\Phi_i$ as above for each of the domain decompositions $\R = (-\infty,-\ell)\cup (-\ell,\infty),  (-\infty, \ell)\cup (\ell,\infty)$. Then $\Phi$ may be obtained as a sum of the $\rho_i(x)\cdot\Psi_i(x,\tau)$, where each $\rho_i$ is a sufficiently localized smooth bump function which is identically one in a neighborhood of a discontinuity of $b$.

\section{Proof of main theorem}\label{s:pfT1}
\subsection{Smooth front profile $u_*$}

In this section, we give the proof where Hypothesis \ref{h:nl} and \ref{h:us} hold.  The proof for Hypothesis \ref{h:up} will follow in the same way with a few alterations and is lined out in Section \ref{ssec:pfdisc}.

Using Lyapunov-Schmidt reduction, we wish to solve $\mc{F}(w,\omega,c) = 0$ for $(\omega,c)$ close to $(\omega_*, c_*)$.  Since $\mL$ has Fredholm index -1, one must alter the setting before the Implicit Function theorem may be applied in the reduction.  In our setting, this alteration is simple.  Let
$$
\mathring{\mX}:= \{ u\in \mX: \la u, \re^{-2\eta\la x\ra}\ra_\mX = 0\}, 
$$  
so that $\mathring{\mX}$ is a closed subset of $\mX$.  It is readily found that $\mL$ maps $\mY$ into $\mathring{\mX}$, so that we may restrict the codomain of our problem.

Furthermore, the linearization $\mL$ has Fredholm index zero when considered as an operator $\mL: \mY\rightarrow \mathring{\mX}$. This follows from Proposition \ref{p:fred} and the fact that, for any exponential weight with $\eta < 0$, the constant function $1$ lies in the cokernel of $\mL$.

Let us define
\beq
\tl\omega = \omega - \omega_*,\qquad \tl c = c - c_*, \quad\Omega= (\tl\omega,\tl c),
\eeq
so that $\mc{F}$ is now a function of $(u; \Omega)\in \mY\times \R^2$. For the following we suppress the dependence on parameters $(\tl \omega,\tl c)$ unless it is needed.

By Hypothesis \ref{h:HC}, for the functions $P_+(x,\tau) := \re^{\ri \tau} p(x), \,\,P_-(x,\tau):= \overline{P_+(x,\tau)}$,
$$
\ker \mL = \textrm{span}\Big\{ P_+,P_-\Big\}.
$$
Furthermore, let us give $u_0\in\ker \mL$ the coordinates 
$$
u_0= a\, P_+ + \overline{a}\,P_-,
$$
with $a,\overline{a}\in \C$.
Next, with the adjoint eigenfunctions $\psi(x)$ and $\overline{\psi(x)}$ as defined in Section \ref{ss:hmr}, we define $\Psi_+(x,\tau):= \re^{\ri \tau} \psi(x), \,\, \Psi_-:=\overline{\Psi_+}$.  The algebraic simplicity assumed in Hypothesis \ref{h:HC} implies
$$
\la P_i, \Psi_j\ra_\mX = \delta_{ij},
$$
where $\delta_{ij}$ is the Kronecker delta.  Then since $\mL$ is Fredholm by Proposition \ref{p:fred}, we have the following decomposition
$$
\mathring{\mX} = \ker \mL^* \oplus \mc{M}, \qquad \mc{M} = \lp( \textrm{span}\{ \Psi_-, \Psi_+\} \rp)^\perp.
$$
This decomposition has the associated projections
$$
\mc{Q}:\mathring{\mX} \rightarrow \ker \mL^*, \qquad \mc{P}:= I - \mc{Q}: \mathring{\mX}\rightarrow \mc{M}.
$$
The projection onto $\ker \mL^*$ can be explicitly defined as
\beq
\mc{Q}u = \sum_{i = \pm} \la u, \Psi_i\ra_\mX\cdot P_i. 
\eeq

Thus, solving $\mc{F}\equiv 0$ is equivalent to solving the following system of equations 
\begin{align}
0 &=  \mc{P}\mF(u_0+u_h,\Omega), \label{e:LS1}\\
0 &= \mc{Q}\mF(u_0+u_h; \Omega), \label{e:LS2}
\end{align}
where $u_0\in \ker \mL^*$ and $u_h \in \mc{M}$.

Since the linearization $\mL$ of $\mF$ about the trivial state $(0;0,0)$ has Fredholm index zero, the linearization of the first equation with respect to $u_h$ is invertible on $\mc{M}$. Since $\mF$ is a smooth function of $(v,\omega, c)$, the Implicit Function theorem guarantees that there exists a smooth function $\varphi:\ker \mL^*\times \R^2 \rightarrow \mc{M}$ such that $u_0+\varphi(u_0;\Omega)$ solves the "auxiliary" equation \eqref{e:LS1} and $\varphi(0;\Omega) = D_{u_0}\varphi|_{(0,\Omega)} = 0$ for $\Omega$ sufficiently small.

Let us expand $\varphi$ about $u_0 \equiv 0$ using the coordinates defined above,
$$
\varphi(u_0;\Omega) = a^2\re^{2\ri t} \varphi_{+}(x;\Omega) + a\overline{a}\,\varphi_0(x;\Omega) + \overline{a}^2 \re^{-2\ri t}\, \varphi_-(x;\Omega) + \mc{O}(|a|^3).
$$
Inserting this into \eqref{e:LS1}, it is readily found that the functions $\varphi_i$ for  $i = 0,-,+$ must solve the differential equations
\begin{align}\label{e:LS3}
\mL(\re^{2\ri t}\varphi_+ ) &=  \re^{2\ri t}(\p_u^2f(x,u_*(x))p^2 )_{xx}, \nonumber\\
\mL(\varphi_0 ) &=  (\p_u^2f(x,u_*(x))p\overline{p} )_{xx},\nonumber\\
\mL(\re^{-2\ri t}\varphi_- ) &=  \re^{-2\ri t}(\p_u^2f(x,u_*(x))\overline{p}^2 )_{xx}.
\end{align}
Note that each solution $\varphi_i$ exists because the right hand side of each equation in \eqref{e:LS3} is exponentially localized in $x$ so that it is an element of $\mathring{\mX}$ and, by the Fredholm alternative, lies in the range of $\mL$.  

Inserting these solutions into the bifurcation equation \eqref{e:LS2}, we obtain the equivalent reduced system of equations
\beq
0 = \Phi_i(u_0; \Omega) := \la \Psi_i,\mF(u_0+\varphi(u_0,\Omega); \Omega) \ra, \quad i = +,-.
\eeq
Then, after calculations similar to \cite[\S VIII.3]{golubitsky88}, the expansion of each $\Phi_i$ about $(a,\overline{a};\tl\omega,\tl c) = (0,0,0,0)$ is found to be
\begin{align}\label{e:bif}
\Phi_+(a,\overline{a}; \Omega) &=   \lp(\lambda_{\tl c}(0) \tl c + \ri\, \tl\omega \rp) a+ \,\theta_+(0,0) \, a|a|^2 +\mc{O}\lp((|\Omega| + |a|^2)|a\|\Omega|+|a|^4\rp),\\
\Phi_-(a,\overline{a}; \Omega) &=  \lp(\overline{\lambda_{\tl c}}(0) \,\tl c - \ri \tl\omega\rp)\overline{a}+ \,\,\theta_-(0,0)\,\,\overline{a}|a|^2+\mc{O}\lp((|\Omega| + |a|^2)|a\|\Omega|+|a|^4\rp),
\end{align}
where $\lambda_{\tl c}(0) = \fr{d \lambda}{d \tl c}(\tl c = 0)\neq 0$ by assumption and
\begin{align}
\theta_+(\tl\omega,\tl c) &= \Big\la \lp(3\p_u^3f(x,u_*) p^2 \overline{p} + \p_u^2f(x,u_*) \lp[p \varphi_0 + \overline{p} \varphi_{+}\rp]\rp)_{xx}, \psi\Big\ra_{L^2_\eta(\R)},\label{e:bif1}\\
\theta_-(\tl\omega,\tl c) &= \Big\la \lp(3\p_u^3f(x,u_*) p \overline{p}^2 +  \p_u^2f(x,u_*) \lp[p \varphi_- + \overline{p} \varphi_{0}\rp]\rp)_{xx}, \overline{\psi}\Big\ra_{L^2_\eta(\R)}.\label{e:bif2} 
\end{align}
As common in Hopf bifurcation, the time-shift symmetry induces complex rotation equivariance of the bifurcation equation, so that we can factor $a$ and obtain an equation that only depends on $|a|^2$. 

Since $\rre\,\lambda_{\tl c}(0)$ is non-zero by assumption, the Implicit Function theorem gives that there exists a bifurcating branch of  solutions $(\tl \omega,\tl c)(|a|^2)$ parameterized by the amplitude $r^2:=|a|^2$ of the coordinate in $\ker \mL^*$. Whenever $\rre\{\theta_+(0,0)\}\neq 0$, one can solve for $\tl c$ as a function of $r$ and readily confirm the statement on the direction of branching in the theorem. 

\begin{Remark}\label{r:qterms}
For the standard Cahn-Hilliard nonlinearity $f(u) = u - u^3$ linearized about $u_* \equiv 0$, we have $f''(u_*) = 0$ so that the expressions for $\theta_{\pm}$ simplify and, in practice, the inhomogeneous problems in \eqref{e:LS3} need not be solved.
\end{Remark}

\begin{Remark}\label{r:par}
We note that instead of restricting the codomain so that $\mL$ is Fredholm index 0, one could also follow the work of \cite{pogan2010} and \cite{sandstede08}  by adding an extra parameter via the ansatz 
$$
u(x) = b \chi(x) + w(x),
$$
where $b\in \R$, $\chi(x) = (1 + \tanh(x))/2$ and $w\in \mX$. When considered in these coordinates, $\mc{F}: \mY\times \R^3\rightarrow \mX$ will then have a linearization which has Fredholm index 0 so that we may then perform a Lyapunov-Schmidt reduction as above. Such an approach would also be necessary when mass conservation only determines the asymptotic mass difference implicitly, say, when mass deposition through the trigger depends on concentrations $\chi=\chi(x,u)$. 
\end{Remark}

\subsection{Alterations for Hypothesis \ref{h:up}}\label{ssec:pfdisc}
The proof under Hypothesis \ref{h:up} follows in a similar manner and we only note the few differences. By taking into account the jump conditions at the discontinuities, it is readily found that $\mF$ and hence $\mL$ maps $\mY^{\#\#}$ into $\mathring{\mX}$.  A routine calculation then shows that the weak derivatives in the right sides of the three equations in \eqref{e:LS3} are well defined and in $\mX$.  The solvability of these equations then follows from the Fredholm Alternative, the fact that the eigenfunction $p$ is exponentially localized, and $\ker (L - \ri k \omega )= \{0\}$ for $k\neq \pm1$.

\section{Instability plateaus --- an explicit example}\label{s:example}

In this section, we study an example where we can establish existence of modulated traveling waves. That is, we are able to verify the assumptions of Theorem \ref{t:hbex}. We first motivate our specific choice of nonlinearity in Section \ref{ss:41}, and then introduce general concepts on absolute and convective instability in bounded domains in Section \ref{ss:42}. Sections \ref{ss:43} to \ref{ss:45} then establish precise asymptotics for the first Hopf instability for long plateaus. Finally, Section \ref{ss:46} concludes by determining the cubic Hopf coefficient and the direction of branching.  

\subsection{Motivation}\label{ss:41}
Establishing the existence of a Hopf bifurcation can generally be cumbersome. While the Hopf bifurcation that we analyzed here is ubiquitous in numerical simulations and experimental observations, it is generally difficult to rigorously prove that the assumptions of our theorem are satisfied. Intuitively, one expects a Hopf instability since for slow speed, mass deposition is slow so that the system develops a long, slowly varying plateau-like state in the intermediate spinodal regime. On this state, one expects a spinodal decomposition instability, with a typical selected spatial wavenumber. Since this instability is stationary in the steady frame, one would expect oscillations in the co-moving frame of the trigger front. The absence of an explicit expression for the trigger front and the lack of tools to detect Hopf eigenvalues makes this problem in general quite intractable. We therefore set up a toy system, where the front is trivial, the ``plateau'' is an actual constant state, and nonlinearities are piecewise constant. As a benefit, we show how to make the above intuition rigorous in terms of branch points and absolute spectra, in particular obtaining corrections to the simple wavenumber prediction from fastest growing modes.

For the remainder of the section we let $u_*(x)\equiv 0$ and study \eqref{e:CH-mf} with nonlinearities of the form $f(x,u) =\chi(x) u + \gamma u^3-\beta u^5$ where $\beta$ and $\gamma$ are real constants with $\beta>0$ and
\beq\label{e:chi}
\chi(x) = 
\begin{cases}
\chi_+ = 1 &\quad x\in [-\ell,\ell]\\
\chi_- = -1&\quad x\in (-\infty, -\ell)\cap (\ell,\infty),
\end{cases}
\eeq
is a triggering mechanism which makes the homogeneous state $u\equiv 0$ linearly unstable inside the interval $[-\ell,\ell]$ and linearly stable everywhere else.  Such triggers have been used to numerically study directional quenching (see \cite{Foard12} and \cite{krekhov}) and are a caricature of many others used in different situations; see Section 1.1 above. Also, by scaling we may assume that $\beta = 1$. This nonlinearity obviously satisfies Hypothesis \ref{h:up} as noted above. Denote
$$
L u =  -\lp(u_{xx} +\chi u\rp)_{xx} + cu_x,
$$
and, as they will be of use in the following propositions, define the constant coefficient operators
\beq\label{e:cstcoop}
L_{\pm}u:= -\lp(u_{xx} +(\chi_{\pm}) u\rp)_{xx} + cu_x .
\eeq

\subsection{Absolute and convective instabilities in bounded domains}\label{ss:42}
One can think of the linear problem with piecewise constant trigger $\chi$ as in (\ref{e:chi}) as a problem on $x\in(-l,l)$ with ``effective'' boundary conditions at $\pm l$, induced by the stable system on either side of the plateau. On the plateau, we see an instability which is advected by the drift term $c\partial_x$ induced by the co-moving frame. Only for sufficiently strong instabilities will the exponential growth outpace the linear advection. The eigenvalue problem on a finite domain is of course ``explicitly solvable'', in principle. On the other hand, calculations very quickly become quite impenetrable and we pursue a more conceptual approach.

In fact, the results in  \cite{Scheel00} and \cite{gluing}  show how to generally compute asymptotic behavior of spectra in finite bounded domains. For large domain length and separated boundary conditions which satisfy a certain non-degeneracy condition, all but finitely many eigenvalues are approximated by a set of curves called the \emph{absolute spectrum}.  This set is determined via the dispersion relation $d(\lambda,\nu)$ obtained by inserting $\re^{\lambda t + \nu x}$ into the asymptotic linearized equation.  By viewing the \emph{temporal} eigenvalue $\lambda$ as a parameter and solving for the \emph{spatial} eigenvalues $\nu = \nu(\lambda)$ ordered by real part $\rre \, \nu_j \geq \,\rre\, \nu_{j+1}$, one finds for well-posed operators that the system has fixed morse index $i_\infty$ so that $\rre \,\nu_{i_\infty}(\lambda) > 0 > \rre \,\nu_{i_\infty+1}(\lambda)$ for all $\lambda$ with large real part.  The absolute spectrum is then defined as 
$$
\Sigma_\rab = \{ \lambda\in \C : \rre\, \nu_{i_\infty}(\lambda) = \rre\, \nu_{i_\infty+1}(\lambda) \}.
$$

Though $\Sigma_\rab$ is not part of the spectrum of the linearized operator on an infinite domain, it dictates whether instabilities saturate the domain or are convected away. For typical problems, $\Sigma_\rab$ has an element $\lambda_\br$ with largest real part which determines when such instabilities arise. It is often the case that $\lambda_\br$ is a branch point of the dispersion relation and hence is an endpoint of a curve in $\Sigma_\rab$ which satisfies $\nu_{i_\infty}(\lambda_\br) = \nu_{i_\infty+1}(\lambda_\br)$; see \cite{rademacher07} or \cite{Scheel00}). Additionally, the results of \cite{Scheel00} give that eigenvalues of the finite domain problem of length $\ell$ accumulate on $\lambda_\br$ with rate $\mc{O}(\ell^{-2})$.  

The work of \cite{gluing} uses these concepts to study the spectrum of a pulse $p(x)$ connecting a stable rest state $p_0$ at $x\rightarrow\pm \infty$ to a plateau state which is close to an unstable rest state $p_1$ for $x\in [-\ell,\ell]$. By viewing such a pulse as the gluing of "front" and "back" solutions between $p_0$ and $p_1$, the limiting spectral set (as $\ell \rightarrow \infty$) of the linearization about this pulse can be decomposed into three parts: the absolute spectrum of the linearization about the unstable state $p_1$, the essential spectrum of the linearization about the state $p_0$, and a finite number of isolated eigenvalues determined by the spectrum of the front and back solutions. Using arguments as in \cite{Scheel00}, it is also shown that an infinite number of eigenvalues converge to the absolute spectrum with $\mc{O}(\ell^{-2})$ rate.

In our setting, the solution $u_*(x)$ can be viewed as a pulse whose asymptotic operator, defined above as $L_-$, has marginally stable spectrum.  We will show for large $\ell$ that eigenvalue crossings are approximated by intersections of the absolute spectrum of $L_+$ with the imaginary axis.  The absolute spectrum, which we denote as $\Sigma_\rab^+$, is determined by the dispersion relation $d_+$ in \eqref{e:dsp1} below.  Furthermore, the first crossing is approximated by where the right-most part of $\Sigma_\rab^+$, which consists of two complex conjugate branch points, intersects $\ri \R$. As the front speed $c$ is decreased, $\Sigma_\rab^+$ moves to the right towards the right half of the complex plane $\C^+$. As discussed above, when $\Sigma_\rab^+\cap \C^+\neq \varnothing$ instabilities which are stronger than the convective motion arise in the domain $[-\ell,\ell]$.  This heuristically indicates that unstable eigenvalues will lie close to $\Sigma_\rab^+\cap \C^+$. In the following, proof of these facts in our specific context is done by hand as the aforementioned results are not directly applicable and do not give explicit expansions of eigenvalues near the branch point.

\subsection{Extended point spectrum}\label{ss:43}

We now begin to verify the spectral hypotheses for our explicit example. In this section, we show that no eigenvalues arise from the front or back solutions.  The genericity of the absolute spectrum, discussed in Section \ref{ss:44}, will then allow us to show that eigenvalues which accumulate onto the absolute spectrum are the first to bifurcate.

In our case, the front and back solutions are $u_*(x)\equiv 0$ which solve the toy problem with $\chi(x)$ defined respectively as
\beq
\chi_f(x)=\begin{cases}
\chi_+  & \quad x\in(-\infty,0]\\
 \chi_- & \quad x\in (0,\infty)
\end{cases}
,\qquad
\chi_b(x) = 
\begin{cases}
\chi_- & \quad x\in(-\infty,0]\\
\chi_+ & \quad x\in (0,\infty)
\end{cases}.
\eeq
These solutions then give the following piecewise-constant coefficient linearizations composed of $L_\pm$
\begin{align}
L_{\mathrm{f}/\mathrm{b}} u: = -\lp(u_{xx} +(\chi_{\mathrm{f}/\mathrm{b}}) u\rp)_{xx} + cu_x, \qquad 
\end{align} 
which have domain $\mY^\#\subset \mX$ where $x_0 = 0$ and $b= \chi_{\mathrm{f}/\mathrm{b}}$.

We analyze the corresponding Evans functions $D_{\mathrm{f}/\mathrm{b}}(\lambda)$ whose zeros are the eigenvalues of $L_{\mathrm{f}/\mathrm{b}}$; for more background see \cite{kapitula13} and references therein.  In this simple case, the Evans function can be expressed in terms of the stable and unstable eigenspaces $E_{\rs}^\pm(\lambda)$ and $E_{\ru}^\pm(\lambda)$ of the first order systems associated with the operators $L_\pm - \lambda$ as in \eqref{e:cstcoop} above. Namely, $D_{\mathrm{f}/\mathrm{b}}(\lambda):= E_{\rs}^\pm(\lambda) \wedge E_{\ru}^\mp(\lambda)$.  Instead of the usual formulation in terms of $u$ and its derivatives, we use a different set of variables in which the jump conditions $(\#\#)$ at $x = \pm \ell$ become continuity conditions.  Namely we let $v = u_x, \theta = u_{xx}+ \chi_{\mathrm{f}/\mathrm{b}} u,$ and $w = \theta_x$ so that the first order systems take the form
\begin{align}\label{e:1st}
u_x &=v \notag\\
v_x &= \theta- \chi_{\pm} u\notag\\
\theta_x&= w\notag\\
w_x&= cv - \lambda u.
\end{align}
The eigenvalues of this system, denoted as $\nu_i^\pm(\lambda)$, are roots of the dispersion relations
\beq\label{e:dsp1}
d_{\pm}(\lambda,\nu) = -\nu^4 -\chi_{\pm}\nu^2 +c \nu - \lambda.
\eeq
We order these roots by decreasing real part 
$$
\rre\{\nu_j^\pm(\lambda)\} \geq \rre\{\nu_{j+1}^\pm(\lambda)\},
$$
and let $$e_{i}^\pm(\lambda):= \lp(1,\,\nu_{i}^\pm(\lambda),\, \nu_{i}^\pm(\lambda)^2+\chi_\pm,\, \nu_i^\pm(\lambda)(\nu_{i}^\pm(\lambda)^2+\chi_\pm)\rp)^T$$ be the corresponding eigenvectors. As mentioned above, for $\lambda$ with large positive real part it can readily be found that
\beq
\rre\{\nu_1^\pm(\lambda)\}\geq \rre\{\nu_2^\pm(\lambda)\} >0 > \rre\{\nu_3^\pm(\lambda)\}\geq \rre\{\nu_4^\pm(\lambda)\}.
\eeq
In fact this splitting holds for all $\lambda$ to the right of $\Sigma_\mathrm{ess}^\pm$, the essential spectrum of $L_\pm$.   For either $j   = 1, 3$, if $\nu_j^i \neq \nu_{j+1}^i$ then $e_j^i$ and $e_{j+1}^i$ span the unstable and stable eigenspaces of  \eqref{e:1st} respectively.  We find up to a normalization factor, for all $\lambda\in \C\diagdown \Sigma_\rab$ with $e_1^\pm \neq e_2^\pm$ and $e_3^\pm \neq e_4^\pm$,
\beq
D_{f}(\lambda) = \det\Bigg| e_{1}^+\, e_{2}^+\, e_{3}^-\,  e_{4}^- \Bigg|,\qquad 
D_{b}(\lambda) = \det\Bigg|  e_{1}^- \,  e_{2}^-\, e_{3}^+\, e_{4}^+ \Bigg|,\label{e:evfcn}
\eeq
where we have suppressed the dependence on $\lambda$ of $e_j^\pm$ inside the determinant.

If for example $\nu_1^- = \nu_2^-$ for some $\lambda_0$, then one must view the spatial eigenvalues as functions of a variable $\zeta$ on a Riemann surface, $\lambda = g(\zeta)$, with a branch point at $\lambda_0$; see \cite[\S 9.1]{kapitula13}.  Since $\nu_1^-$ is analytic in $\zeta$, the vectors $e_1^-$ and $\fr{d}{d\zeta}e_1^-$ form a basis for the corresponding unstable eigenspace. 

With these definitions, we readily obtain the following lemma.

\begin{Lemma}\label{l:epsp}
For all speeds $c>0$, the functions $D_f(\lambda)$ and $D_b(\lambda)$ have no zeros in the set $\C \diagdown \Sigma_\rab^+$.  Furthermore, they have a non-vanishing limit as $\lambda$ approaches $\Sigma_\rab^+$.
\end{Lemma}
\begin{Proof} As the argument will be the same for the back, we only consider the front. By applying Sobolev embeddings to the numerical range of both $L_\pm$, taking care to mind the jump conditions $(\#)$, it is readily found that $L_f$ is uniformly sectorial on $L^2(\R)$ in the plateau length $\ell$.  This implies that both $D_f$, being analytic off of the absolute spectrum, does not vanish identically in any connected component of $\C \diagdown\Sigma_\rab$.

It can be readily found that $D_f(0)\neq 0$. Assuming that $\lambda\neq0$, we split the proof of the first statement into two cases.

\textbf{Case 1:} Assume that $\lambda$ is such that $\nu_1^+ \neq \nu_2^+$  and $\nu_3^-\neq \nu_4^-$.

In this case the Evans functions $D_{f/b}$ are given by \eqref{e:evfcn} above,
\begin{align}
D_{f}(\lambda) &= \det\left(\begin{array}{cccc}1 & 1 & 1 & 1 \\\nu_{1}^- & \nu_{2}^-  & \nu_{3}^+ & \nu_{4}^+ \\(\nu_{1}^-)^2+\chi_- & (\nu_{2}^-)^2+\chi_-  & (\nu_{3}^+)^2+\chi_+  & (\nu_{4}^+)^2+\chi_+  \\
\nu_{1}^-( (\nu_{1}^-)^2+\chi_-) & \nu_{2}^-( (\nu_{2}^-)^2+\chi_-) & \nu_{3}^+( (\nu_{3}^+)^2+\chi_+) & \nu_{4}^+( (\nu_{4}^+)^2+\chi_+)\end{array}\right)\notag\\
&= \det \left(\begin{array}{cccc}1 & 1 & 1 & 1 \\\nu_{1}^- & \nu_{2}^- & \nu_{3}^+ & \nu_{4}^+ \\  (\nu_{1}^- )^2 & (\nu_{2}^- )^2 &  (\nu_{3}^+)^2 &  (\nu_{4}^+)^2 \\ (\nu_{1}^- )^3 & (\nu_{2}^- )^3 &  (\nu_{3}^+)^3 &  (\nu_{4}^+)^3\end{array}\right),
\end{align}
a Vandermonde determinant which we shall denote as $V(v_1^-,v_2^-,v_3^+,v_4^+)$. This equality can be obtained using the dispersion relation to find $(\nu_i^\pm)^2 + \chi_\pm = -\fr{\lambda - c\nu_i^\pm}{(\nu_i^\pm)^2}$ and then performing elementary row operations.

Hence, $D_f(\lambda) = 0$ if and only if $\nu_{3}^-(\lambda) = \nu_{2}^+(\lambda)$.  This means that both dispersion relations $d_\pm$ are satisfied simultaneously and, since $\chi_+ \neq \chi_-$, that $\nu_{3}^- = \nu_{2}^+ = 0$. But we also have that $\nu_i^\pm(\lambda) = 0$ if and only if $\lambda = 0$.  Therefore $D_f(\lambda)\neq0$.

\textbf{Case 2:} Assume either $\nu_1^- = \nu_2^-$ or $\nu_3^+ = \nu_4^+$. 

Say only the latter holds.  Then we have
$$
D_f(\lambda) = \p_\zeta \nu_3^+ \cdot \fr{d}{d\nu_4^+}\,\Big|_{\nu_4^+ = \nu_3^+} V \neq 0,
$$
for all $\lambda\neq0$ because once again $D_f(\lambda) = 0$ if and only if $\nu_3^+ = \nu_2^-$ which holds if and only if $\lambda = 0$. If $\nu_1^- = \nu_2^-$ then take the derivative of $V$ with respect to $\nu_1^-$. If both equalities hold then take the derivatives of $V$ with respect to both $\nu_1^-$ and $\nu_3^+$. Note, we have that $\p_\zeta \nu_3^+\neq 0$ because double roots are simple for all $c>0$.  This gives the proof of the first part of the lemma.

To prove the second statement we note that if $\lambda \rightarrow \lambda_0\in\Sigma_\rab^+$ then, by definition, $\rre\,\nu_2^+- \rre\,\nu_3^+ \rightarrow 0$.  If $D_f(\lambda)$ were to approach zero as well, then arguments used above give that $\lambda_0 $ must be 0, which is readily found to lie in the complement of the absolute spectrum for all speeds $c>0$. This gives the proof of the second statement and completes the lemma.

\end{Proof}

\subsection{Branch points, rescalings and asymptotics}\label{ss:44}
We now analyze the dispersion relation near the rightmost point of the absolute spectrum in more detail. We give explicit formulas describing how branch points cross the imaginary axis and how the spatial eigenvalues $\nu(\lambda)$ behave around them.  Furthermore, we will study how $\Sigma_\rab^+$ behaves near $\lambda_\br(c)$.

For $c = c_\rlin$, it is readily found that the essential spectrum, $\Sigma_\mathrm{ess}^+$, lies in the closed left-half plane when considered in an exponentially weighted space with weight $\re^{\mu_\rlin x}$.  Since $\Sigma_\rab^+$ generically lies to the left of $\Sigma_\mathrm{ess}^+$, elementary calculation shows that, for $c$ near $c_\rlin$, the right most part of the absolute spectrum consists of a pair of complex conjugate branch points of the dispersion relation $d_+$. Such branch points, which we denote as $ \lambda_\br(c), \overline{\lambda_\br(c)}$, solve the algebraic system
 \begin{align}
d_+(\lambda,\nu) &= 0,\\
\fr{\p}{\p\nu} d_+(\lambda,\nu) &=0,
\end{align}
for some double spatial eigenvalue which we denote as $\nu_\br(c):= \nu(\lambda_\br(c))$.  

In the context of front invasion into an unstable state, if $\nu_\br(c)$ satisfies what is known as a "pinching"condition, the speed $c = c_\rlin$ for which $\lambda_\br(c)\in \ri \R$ is called the \emph{linear spreading speed}; see  \cite{brevdo96}, and \cite{holzer14}. Such "pinched double root" solutions of the Cahn-Hilliard dispersion relations have been studied previously and explicit expressions for $\lambda_\rlin := \lambda_\br(c_\rlin)$ and $\nu_\rlin := \nu_\br(c_\rlin)$ have been obtained. As they will be of use in the following, we sum them up in the following lemma.
\begin{Lemma}\label{l:calc}Given $f$ and $u_*$ as above, for $\alpha = \p_uf'(0,u_*(0))$, we have the following 
\begin{align}
\lambda_\rlin &= \ri (3+ \sqrt{7})\sqrt{\fr{2+\sqrt{7}}{96}}\cdot \alpha^2\notag\\
c_\rlin &= \fr{2}{3\sqrt{6}}(2+\sqrt{7}) \sqrt{\sqrt{7}-1} \cdot \alpha^{3/2}\notag\\
\mu_\rlin&:= \rre\{ \nu_\rlin\} = - \sqrt{\fr{\sqrt{7} - 1}{24}} \cdot \alpha^{1/2}\notag\\
\kappa_\rlin&:= \rim\{\nu_\rlin\} = \sqrt{\fr{\sqrt{7} + 3}{8}}\cdot \alpha^{1/2}.
\end{align}
\end{Lemma}
\begin{Proof}
These quantities can be found in \cite[Lem 1.3]{Scheel12} or \cite{vanSaarloos03}. 
\end{Proof}

In the next section, we will use spatial dynamics to obtain precise expansions for the first eigenvalue crossing and its corresponding eigenfunction. In order to do this we must obtain expansions for the spatial eigenvalues which solve the dispersion relation \eqref{e:dsp1} for $\lambda$ near $\lambda_\rlin$. Thus let $\hat \lambda = \lambda - \lambda_\rlin$, $\hat\nu = \nu - \nu_\rlin$, $\hat c = c - c_\rlin$ and $\Lambda = (\hat\lambda, \hat c)$.  In these variables the dispersion relation \eqref{e:cjeqn} takes the form
\beq\label{e:cjdr}
\hat d_+(\hat\lambda, \hat \nu) := \hat \nu^4+4\nu_\rlin\hat \nu^3+(1+6\nu_\rlin^2)\hat\nu^2   - \hat c \hat \nu + \hat\lambda - \hat c \nu_\rlin.
\eeq
We characterize the roots $\hat\nu(\hat\lambda, \hat c)$ in the following lemma.
\begin{Lemma}\label{l:cjevals}
 The dispersion relation \eqref{e:cjdr} has four roots, $\hat\nu_{\rs},\hat\nu_{\ru},\hat\nu_{\rcs},\hat\nu_{\rcu}$, which are functions of $\Lambda\in \C\times \R$ and, for all $\Lambda$ close to $(0,0)$, satisfy the following properties
 \begin{enumerate}
 \item $\hat\nu_{\rs/\ru} = -2\nu_\rlin \pm  \sqrt{-2\nu_\rlin^2 - 1} + \mc{O}(|\Lambda|).$
 
 \item The roots $\hat\nu_{\rcs/\rcu}$ solve
\beq\label{e:redpol}
 \hat \nu^2 + b_1(\Lambda) \hat\nu + b_0(\Lambda) = 0, \quad 
\eeq
where, setting $\gamma_\rlin = (1+6\nu_\rlin^2)$, the coefficients $b_0$ and $b_1$ are analytic functions of $\Lambda$ with leading order expansions
\beq
b_1(\Lambda) =\lp( \fr{4\nu_\rlin}{\gamma_\rlin^2} - \fr{1}{\gamma_\rlin}  \rp) \hat c  - \fr{4\nu_\rlin}{\gamma_\rlin^2}\hat\lambda+ \mc{O}(|\Lambda|^2),\qquad  b_0(\Lambda) = -\fr{\nu_\rlin}{\gamma_\rlin}\hat c + \fr{1}{\gamma_\rlin} \hat \lambda + \mc{O}(|\Lambda|^2).
\eeq

\item  For all $\Lambda$ with $\hat\lambda +\lambda_\rlin \not \in \Sigma_\rab^+$ the roots $\hat\nu_{\rcs/\rcu}$ split in the following way
\beq\label{e:nuspl}
\rre\{\hat \nu_\rcs\} < -\fr{b_1(\Lambda)}{2} < \rre\{\hat\nu_\rcu\}.
\eeq

 \end{enumerate}
\end{Lemma}

\begin{Proof}
Property (i) is easily proved using standard perturbation techniques.  Property (ii) is obtained using multi-parameter expansions and the Weierstrass Preparation Theorem; see for example \cite[Ch. 4]{seyranian03}.  We note that \eqref{e:redpol} may be used to determine the branch point  $(\hat\lambda_\br(\hat c), \hat\nu_\br(\hat c))$ in the shifted dispersion relation \eqref{e:cjdr},  for $\hat c$ near zero.  Indeed, $\hat\lambda_\br(\hat c)$ must satisfy
\beq\label{e:nuh}
0=b_0(\hat\lambda_\br(c),\hat c) - \fr{b_1(\hat\lambda_\br(\hat c),\hat c)^2}{4}, 
\eeq
and hence has the form $\hat\lambda_\br(\hat c) = \nu_\rlin \hat c + \mc{O}(\hat c^2),$
while $\hat\nu_\br(c) = -\fr{b_1(\lambda_\br(\hat c),\hat c)}{2}$.

It now remains to prove property (iii).  In order to find expansions for the roots of \eqref{e:redpol}, we make the change of variables $\hat\nu = \tl \nu - \fr{b_1(\Lambda)}{2}$ so that 
\beq\label{e:redpol2}
0 = \tl\nu^2 + \beta(\Lambda), \qquad\text{with}\quad \beta(\Lambda) = -b_0(\Lambda)+b_1(\Lambda)^2.
\eeq
Fixing $\hat c$, setting $\tl \lambda = \hat \lambda - \hat\lambda_\br(\hat c)$,  and expanding near $\hat\lambda_\br(\hat c)$ we obtain 
\beq
0= \tl \nu^2 + \tl \lambda\,\tl b_2(\tl\lambda, \hat c),
\eeq
for some function $\tl b_2$ which is analytic in $\tl\lambda$ with $\tl b_2(0,0) = \fr{1}{\gamma_\rlin}$. Finally, setting $\tl \lambda = -\zeta^2$ and scaling $\tl\nu_1 = \tl \nu \zeta$ we find 
\beq
\tl \nu_1 = \pm \sqrt{\tl b_2(-\zeta^2,\hat c)} = \pm \gamma_\rlin^{-1/2} + \mc{O}(|\zeta^2| + |\hat c|).
\eeq
Unwinding all of these scalings gives two roots, $\hat \nu_{\rcu}$ and $\hat\nu_\rcs$, which are analytic on the Riemann surface defined by $\zeta$, and satisfy
$$
\rre\{\hat \nu_\rcs\} < -\fr{b_1(-\zeta^2,\hat c)}{2} < \rre\{\hat\nu_\rcu\},\quad\text{for all}\quad\zeta \not\in S_\rab = \{\xi: \xi^2\,\tl b_2(-\xi^2, \hat c) \in\R_-\}, 
$$
where $\R_-$ is the non-positive part of the real line.  This completes the proof of the lemma.
\end{Proof}

We remark that the calculations of Lemma \ref{l:calc} imply that  for all $\Lambda$ small, the eigenvalues $\hat\nu_{\rs/\ru}$ are bounded away from the imaginary axis, with real parts of opposite sign.

The following lemma shows that $\Sigma_\rab^+$ is generic near the branch point $\lambda_\br(c)$ for all $c$ near $c_\rlin$. 
The result of this lemma is the reducibility hypothesis in \cite[\S7 ]{gluing}. Coupled with Lemma \ref{l:epsp}, this will imply that bifurcating spectra of $L$ are only found near $\Sigma_\rab^+$.

\begin{Lemma}\label{l:gen}
Let $V\subset (\C  - \Sigma_\mathrm{ess})$ be an open, bounded, and connected set containing the branch point $\lambda_\mathrm{br}(c)$ for all $c$ close to $c_\rlin$.  Given such a speed $c$, each $\lambda\in (\Sigma_\rab^+\cap V)\diagdown \{\lambda_\br(c)\}$ satisfies the following:
\beq
\nu_{\ri_\infty}(\lambda)\neq \nu_{\ri_\infty+1}(\lambda),\qquad \fr{d(\nu_{\ri_\infty} -  \nu_{\ri_\infty+1})}{d\lambda}\neq 0.
\eeq
\end{Lemma}
\begin{Proof}

 By definition, for any $\lambda \in \Sigma_\rab\cap V$ there exist $\nu\in \C$ and $\gamma \in \R$ such that 
$$
d_+(\lambda, \nu) = d_+(\lambda,\nu+\ri\gamma) = 0.
$$
Expanding from the branch point we find, after the change of variables $(\tl \lambda,\tl \nu) = (\lambda - \lambda_\br, \nu - \nu_\br)$, that $\tl \lambda + \lambda_\br\in \Sigma_\rab^+$ satisfies
\begin{align}\label{e:abgen}
\tl \lambda &= b \tl\nu^2 + \mc{O}(\tl\lambda\tl\nu, \tl \lambda^2, \tl \nu^2),\\
 \tl \lambda &= b(\tl\nu + \ri \gamma)^2 + \mc{O}(\tl\lambda\tl\nu, \tl \lambda^2, \tl \nu^2),
\end{align}
where $b\in \C$ is a non-zero constant.
This implies that 
\beq
\tl\nu = -\fr{\gamma}{2} \ri + \mc{O}(\gamma^2).
\eeq
By substituting this into the first equation of \eqref{e:abgen} we then find
\beq
\tl\lambda = -\fr{\gamma^2}{4} + \mc{O}(\gamma^3).
\eeq
which implies for $0<\tl \lambda <<1$ that $\gamma \neq 0$ and 
\beq
\fr{d(\tl \nu_{\ri_\infty} - \tl \nu_{\ri_\infty+1})}{d\tl\lambda}\neq 0,
\eeq
where $\ri_\infty$ denotes the Morse index of the first order system corresponding to $L_+$, and counts the dimension of the unstable eigenspace as $\lambda \rightarrow \infty$.
\end{Proof}

\subsection{Spatial dynamics near the branch point}\label{ss:45}
Having collected spectral facts in Sections \ref{ss:42} - \ref{ss:44}, we now are able to use spatial dynamics to characterize the first eigenvalue crossing and its corresponding eigenfunction.  We construct eigenfunctions of $L$ by conjugating with $\re^{\nu_\rlin x}$ and solving the finite domain eigenvalue problem for $x\in [-\ell,\ell]$ subject to boundary conditions induced by the dynamics for $x\in \R \diagdown [-\ell,\ell]$.  

Inserting $u = \re^{\nu_\rlin x} \tl u$ into the eigenvalue equation $L u  = \lambda u$, dividing by $\re^{\nu_\rlin x}$, and using the fact that $d_+(\lambda_\rlin, \nu_\rlin) = \fr{d}{d\nu} d_+(\lambda_\rlin, \nu_\rlin) = 0$,  we obtain an equivalent eigenvalue equation which, when expressed in scaled variables, takes the form
\begin{align}
\p_x^4\tl u + 4\nu_\rlin \p_x^3 \tl u+(\chi_++6\nu_\rlin^2)\p_x^2\tl u - \hat c \p_x\tl u  +(\hat \lambda - \hat c \nu_\rlin) \tl u &=0,\quad x\in [-\ell,\ell],\label{e:cjeqn}\\
\p_x^4\tl u + 4\nu_\rlin \p_x^3\tl u +(\chi_-+6\nu_\rlin^2)\p_x^2 \tl u - (\hat c+ 2(\delta \chi)\nu_\rlin)\p_x \tl u+(\hat\lambda - \hat c\nu_\rlin - (\delta\chi)\nu_\rlin^2)\tl u&=0,\quad x\in \R \diagdown [-\ell,\ell],\label{e:cjeqn1}
\end{align}
where $\delta\chi = \chi_- - \chi_+.$
Using the coordinates of \eqref{e:1st}, these operators have the first order systems
\begin{align}\label{e:1st1}
\tl  u_x &=\tl v - \nu_\rlin \tl u \notag\\
\tl  v_x &= \tl\theta- \chi_{\pm} \tl u - \nu_\rlin \tl v\notag\\
\tl  \theta_x&= \tl w - \nu_\rlin \tl \theta\notag\\
\tl  w_x&= (c_\rlin+\hat c)\tl v - (\lambda_\rlin+\tl \lambda) \tl  u - \nu_\rlin \tl  w.
\end{align}
If $\hat\nu_i^\pm$ are the eigenvalues for this system, ordered by decreasing real part, then the corresponding eigenvectors take the form
$$
\hat e_i^\pm = \lp(1,\hat \nu_i^\pm + \nu_\rlin ,\, (\hat \nu_i^\pm + \nu_\rlin )^2+\chi_\pm,\, (\hat \nu_i^\pm + \nu_\rlin )\lp( (\hat \nu_i^\pm + \nu_\rlin )^2+\chi_\pm \rp) \rp)^T, \quad i = 1,2,3,4.
$$ 
Note, with $\chi_+$ chosen, the eigenvalues of \eqref{e:1st1} are precisely the scaled spatial eigenvalues derived in Lemma \ref{l:cjevals} above. Also, for $\lambda + \hat \lambda_\rlin \in \C \diagdown \Sigma_\mathrm{ess}^-$, the subspaces $\hat E_-^{\rs}:= \mathrm{span}_{i = 3,4} \{ e^-_{i}\} $ and $\hat E_-^{\ru}:=\mathrm{span}_{i = 1,2} \{ e^-_{i}\}$ are the stable and unstable eigenspaces of \eqref{e:1st1} with $\chi_-$ chosen.
 
The boundary conditions at $x = \pm l$ for the eigenfunction are determined as follows. In order for $\tl u$ to be an $L^2(\R)$ eigenfunction, it is necessary and sufficient to require exponential decay as $|x|\rightarrow \infty$. Hence, for $\tl U:=(\tl  u, \tl  v, \tl  \theta, \tl  w)^T$, we require
 \begin{align}\label{e:expbc}
\tl U(\ell) \in \tl  E_-^{\rs},\qquad \tl  U(- \ell) \in \tl  E_-^{\ru}.
\end{align} 

We note that the dimensions of the boundary spaces $\tl  E_-^{\rs/\ru}$ are the same as the corresponding subspaces for the unconjugated problem $L_- u =  \lambda u$. This can be seen by homotoping the conjugation factor $\re^{\nu_\rlin s \,x}$ from $s = 0$ to $s = 1$ and noticing that the essential spectrum of $L_-$ never intersects some sufficiently small neighborhood of $\lambda_\rlin$, implying that no spatial eigenvalue $\nu_i^-$ crosses the imaginary axis during this homotopy. 

Finally, let $\tl  E_+^\rcs$ be the 2-dimensional eigenspace of \eqref{e:1st1} (with $\chi_+$ chosen) spanned by the eigenvectors of $\hat \nu_\rcs$ and $\hat \nu_\rs$. Define $\tl  E_+^\rcu$ in the same way so that it is spanned by the eigenvectors associated with $\hat\nu_\rcu$ and $\hat \nu_\ru$. We remark that both of these subspaces are analytic in the Riemann surface variable $\zeta$ used in the proof of Lemma \ref{l:cjevals} and can be analytically continued as $\zeta$ approaches $S_\rab$, also defined in the above proof.

With these definitions we obtain the following lemma which precludes embedded eigenvalues (see \cite[\S5.3]{Scheel00}),  and will be important in the construction of eigenfunctions. 
\begin{Lemma}\label{c:nondeg}(Non-Degenerate Boundary Conditions) For all $\Lambda$ close to $(0,0)$ with $\hat\lambda+\lambda_\rlin\not\in \Sigma_\rab$, the conjugated eigenspaces of $L_\pm$ satisfy  
\beq\label{e:tranint}
\tl  E_-^{\ru}\pitchfork \tl  E_+^\rcs  = \{0\},   \qquad \tl  E_-^{\rs} \pitchfork \tl  E_+^\rcu = \{0\}, 
\eeq
where $\pitchfork$ denotes the transverse intersection of linear subspaces.
\end{Lemma}
\begin{Proof}
This follows from Lemma \ref{l:cjevals} using similar arguments as in Lemma \ref{l:epsp}. 
\end{Proof}

We are now able to state our existence result and give expansions for the first crossing eigenvalues and their eigenfunctions.  This is done in the following proposition.  
\begin{Proposition}\label{p:efcns}

For $\ell>0$ sufficiently large, there exists a speed $c_*>0$ and simple eigenvalues $\lambda_*(c,\ell),\overline{\lambda_*(c,\ell)}$ of $L$ with the following properties for $c\sim c_*$:
\begin{enumerate}
\item\emph{(First Crossing)} There exists some $\epsilon>0$ so that for all $c>c_\rlin -\epsilon$, $\lambda_*(c_*,\ell) $ and $\overline{\lambda_*}(c_*,\ell)$ are the only eigenvalues lying in the closed right half-plane. 

\item\emph{(Bifurcation)} $\lambda_*(c,\ell)$ is an algebraically simple eigenvalue and satisfies $$\lambda_*(c_*,\ell)  = \ri \kappa_*(c_*,\ell), \quad \fr{d\rre\{\lambda_*\}}{d c}|_{c = c_*}<0.$$
\item\emph{(Expansions)} For $\hat c \in \R$ and $\hat\lambda \in \ri \R$, the crossing speed $c_*(\ell) = c_\rlin + \hat c$ and crossing location $\lambda_*(c_*,\ell)  = \lambda_\rlin+\hat\lambda$ satisfy
\beq\label{e:efcons}
\hat \lambda =  \ri\fr{\pi^2}{4\mu_\rlin \ell^2}(-1 +6(\mu_\rlin^2+\kappa_\rlin^2))+\mc{O}(\ell^{-3}) , \qquad \hat c = - \fr{\pi^2}{4\mu_\rlin \ell^2 }(1+6(\mu_\rlin^2 - \kappa_\rlin^2))+\mc{O}(\ell^{-3}),
\eeq
with $\kappa_\rlin := \rim\{\nu_\rlin\}$ and $\mu_\rlin := \rre\{\nu_\rlin\}$.
\end{enumerate}

Associated with $\lambda_*$, $L$ has an eigenfunction $p$ and corresponding adjoint eigenfunction $\psi$, which satisfy the following properties:
\begin{enumerate}\setcounter{enumi}{3}
\item For $x\in [-\ell,\ell]$,
\begin{align}\label{e:efcn}
p(x) &= A \re^{(\nu_\rlin+\alpha(\ell)) x}\lp( \sin\lp(\fr{\pi \lp(x - \ell\rp)}{2\ell}\rp)+\mc{O}(\ell^{-1})\rp),\\
  \psi(x) &= B \re^{-(\overline{\nu_\rlin+\alpha(\ell)}) x} \lp(\sin\lp(\fr{\pi (x-\ell)}{2l}\rp)+\mc{O}(\ell^{-1})\rp).
\end{align}
Furthermore, for $j = 1,2,3$
\begin{align}
\p_x^j p(x) &= (\nu_\rlin + \alpha(\ell))^j p(x) + \mc{O}(\ell^{-1}),\\
\p_x^j \psi(x) & = (\nu_\rlin + \alpha(\ell))^j \psi(x) + \mc{O}(\ell^{-1}).
\end{align}
Here the error terms are uniform in $x$, $\alpha(\ell) = \mc{O}(\ell^{-2})$, and $A,B>0$ are undetermined normalization constants.

\item Let $U_h := (h,h_x,\; h_{xx}+\chi_- h, \; h_{xxx}+(\chi_- h)_x)^T$ as in \eqref{e:1st} above. Then for $h = p$ or $h = \psi$, there exists a constant $C>0$, independent of $\ell$ such that,
 \begin{align}
 |U_h(x)| &\leq C \ell^{-1}\re^{-\mu_\rlin\ell}\re^{\delta (x+\ell)}, \quad \;\;x\leq-\ell,\label{e:efcndec}\\
 |U_h(x)| &\leq  C\ell^{-1}\re^{\mu_\rlin\ell}\re^{-\delta' (x-\ell)}, \quad\;\;x\geq \ell,\label{e:efcndec2}
 \end{align}
with $\delta = |\rre\{\nu_2^-(\lambda_*)\}| >0$, $\delta' = |\rre\{\nu_3^-(\lambda_*)\}|>0$, and $\nu_i^-(\lambda)$ defined in \eqref{e:dsp1} above. 

\end{enumerate}
\end{Proposition}

\begin{Proof}

Existence of $\lambda_*$ and properties (ii), (iii), and (iv) will all follow from our construction of a solution to the first order system associated with \eqref{e:cjeqn}.  

Property (i) follows using similar methods as in \cite[\S 6]{gluing} and the fact that $\lambda_\br(c)$ is the right-most part of $\Sigma_\rab^+$ for all $c$ near $c_*$. In particular, for $V$ as in Lemma \ref{l:gen}, possibly enlarged to contain the positive real part of the sector which contains the spectrum of $L$, a construction similar to the following can be used to obtain that any $\lambda\in V$ not in a sufficiently small neighborhood of $\Sigma_\rab^+$ is not an eigenvalue.

Now let us begin our construction of the eigenfunctions on the interval $[-\ell,\ell]$. Since the construction of the adjoint eigenfunction $\psi$ follows in the same way, we only describe how $p$ is obtained. After an analytic change of variables, the first order system \eqref{e:1st1} can be split into hyperbolic and center dynamics as
\begin{align}
\dot W_\mathrm{h} &= \left(\begin{array}{cc}\hat\nu_\rs & 0 \\0 & \hat\nu_\ru\end{array}\right)w_\mathrm{h},\label{e:hsys}\\
\dot W_\rc &= \left(\begin{array}{cc}\alpha & 1 \\-\beta & \alpha\end{array}\right) w_\rc\label{e:csys},
\end{align}
where $W := (W_\mathrm{h},W_\rc)^T\in \C^2\times \C^2$, and $\alpha = \alpha(\Lambda), \beta = \beta(\Lambda)$ are parameters, analytic in $(\hat\lambda,\hat c)$, which unfold the two-dimensional Jordan block at $\Lambda = (0,0)$; see \cite{arnold}, for instance. In fact, the leading order expansion of $\beta$ is given in \eqref{e:redpol2} and $\alpha = -b_1(\Lambda)/2$. This can be seen by comparing the characteristic polynomial of the matrix on the right hand side of \eqref{e:csys} to \eqref{e:redpol}.

Next we study how the non-degeneracy conditions given by Lemma \ref{c:nondeg} affect the boundary conditions in \eqref{e:expbc}. Since $\tl E_{\rcs}^+ \oplus \tl E_{\rcu}^- = \C^4$, there exist linear transformations $T_\pm$ such that, for $W_{\mathrm{h}}(x) = (w_\rs,w_\ru)^T$ and $W_\rc = (w_{\rc,0}, w_{\rc,1})^T$,
\begin{align}\label{e:hypBC0}
\left(\begin{array}{c}w_\rs(-\ell) \\w_{\rc,0}(-\ell)\end{array}\right) = T_- \left(\begin{array}{c}w_\ru(-\ell) \\w_{\rc,1}(-\ell)\end{array}\right),\qquad \left(\begin{array}{c}w_\ru(\ell) \\w_{\rc,1}(\ell)\end{array}\right) = T_+\left(\begin{array}{c}w_\rs(\ell) \\w_{\rc,0}(+\ell)\end{array}\right).
\end{align}

Then, given the flow $\Phi_{x,y}$ of the system \eqref{e:hsys} - \eqref{e:csys}, any solution must satisfy the matching condition $\Phi_{\ell,-\ell} W(-\ell) = W(\ell).$  Using a Lyapunov-Schmidt reduction (i.e. project onto the stable, unstable, and center subspaces) it can be obtained for some constants $c_1,c_2\in\C$ that 
\beq\label{e:hypBC1}
w_\rs(-\ell) = c_1 w_{\rc,0}(-\ell) + \mc{O}(\re^{-\delta \ell}),  \quad  w_\ru(\ell) = c_2 w_{\rc,1}(\ell) + \mc{O}(\re^{-\delta \ell}),
\eeq
from which boundary conditions on $w_\ru(-\ell)$ and $w_\rs(-\ell)$ can be determined via \eqref{e:hypBC0}. The solvability of this reduction follows from Lemma \ref{c:nondeg}. 

Furthermore, we readily obtain
\beq\label{e:hypBC2}
|W_\mathrm{h}(\pm \ell)| \leq C |W_\rc(\pm\ell)|+ \mc{O}(\re^{-\delta \ell}),
\eeq 
for some $\tl\delta>0$ dependent on $\hat\nu_\rs$ and $\hat\nu_\ru$. The boundary conditions for $W_\rc(\pm \ell)$ can then be obtained from this reduction and can be found to be independent of $W_\mathrm{h}(\pm \ell)$ up to an $\mc{O}(\re^{-\tl\delta \ell})$ correction.

We now construct solutions to the center system \eqref{e:csys}. We make the scalings $x = \tl x - \ell$,  $W_\rc = \re^{\alpha \tl x}\tl W_\rc$, and define $\mu^2 = -\beta $ so that \eqref{e:csys} becomes
\beq\label{e:ceqn}
\tl W_\rc' = \left(\begin{array}{cc}0 & 1 \\-\mu^2& 0\end{array}\right) \tl W_\rc.
\eeq
In order to ease the derivation, the boundary conditions on the center system may,  up to an $\mc{O}(\re^{-\delta \ell})$ correction, be written as
$$
\left(\begin{array}{c}-1 \\ r_+\end{array}\right)^T \tl W_\rc(2\ell) = 0,\qquad \left(\begin{array}{c}-1 \\r_-\end{array}\right)^T \tl W_\rc(0) = 0,
$$
for some $r_\pm\in \C$ which depend analytically on $\mu$.

Under these conditions, the system \eqref{e:ceqn} has the solution, $\tl W_{\rc,*}$, with first component $w(\tl x) = A(\sin(\mu x) + r_-\mu \cos(\mu x))$, where $\mu$ must satisfy the equation
\beq
\fr{\tan(2\mu \ell)}{2\mu \ell} = \fr{(r_+ - r_-)}{2\ell(1+r_-\,r_+\mu^2)}.
\eeq
For $l$ large, this equation has the solutions $\mu = \fr{\pi k}{2\ell} + \mc{O}(\ell^{-2})$ for integers $k\neq 0$.  To obtain the first eigenvalue crossing, we set $\mu = \fr{\pi}{2\ell}$ so that 
$$
\tl W_{\rc,*}(\tl x) = A\left(\begin{array}{c}\sin(\fr{\pi \tl x}{2\ell})  \\\fr{\pi}{2\ell}\cos(\fr{\pi \tl x}{2\ell}) \end{array}\right)+\mc{O}(\ell^{-1}),
$$
with error term uniform in $x$.  

Recalling that $\mu^2 = -\beta$, we have
$$
\fr{\pi^2}{4\ell^2} =-\beta(\Lambda) = \fr{1}{\gamma_\rlin}(\hat\lambda - \nu_\rlin \hat c) + \mc{O}(|\Lambda|^2),
$$
which can then be solved for $\hat\lambda\in \ri\R$ and $\hat c\in \R$ to obtain the expressions in \eqref{e:efcons} as desired.  Inserting these expressions into the conjugating exponent $\alpha$, we find it has the asymptotics
$$
\alpha(\ell) = -\fr{b_1(\Lambda)}{2} = \mc{O}(\ell^{-2}).
$$

The eigenfunction $p$ given in the statement of the proposition can then be obtained by unwinding all the scalings made above.

To obtain the decay conditions in (v), we notice that $|\tl W_{\rc,*}(0) |, |\tl W_{\rc,*}(2\ell) | \leq C' l^{-1}$ and thus, given the estimates \eqref{e:hypBC2}, for $h=\psi, p$ and some constant $C'>0$,
$$
|U_h(\pm l)|\leq C' \ell^{-1}\re^{\pm\mu_\rlin\ell}.
$$
The boundary conditions \eqref{e:expbc} then give the estimates \eqref{e:efcndec} and \eqref{e:efcndec2} above.

\end{Proof}

\begin{Remark}\label{r:numer}
The first and second leading order terms in the expansions for the Hopf crossing location $\lambda_*$ and speed $c_*$ in \eqref{e:efcons} were compared with numerical calculations of the spectrum of $L$ and were found to be in excellent agreement.  The operator was considered on a large but finite domain and was discretized using fourth-order accurate finite differences. 
\end{Remark}

\begin{Remark}\label{r:singpert}
We note that the above result could also be obtained via a similar, and in some sense equivalent, geometric singular perturbation method.  If $\hat\lambda$ and $\hat c$ are scaled by $\epsilon = \ell^{-2}$, one obtains a slow equation which is equivalent to a heat equation.  Furthermore for $\epsilon = 0$ the boundary conditions for a solution on the slow manifold reduce to Dirichlet conditions. By solving this system, the same leading order expansions for the eigenfunction, and eigenvalue-front speed pair $(\lambda_*,c_*)$ may be obtained.
\end{Remark}

\subsection{Nonlinear Hopf bifurcation --- direction of branching}\label{ss:46}

We are now able to state our main result of this section which gives the existence of bifurcated solutions and determines the direction of bifurcation in terms of the the cubic nonlinearity parameter $\gamma$.

\begin{Theorem}\label{t:exam}
For $f$ and $u_*$ described above and $\ell>0$ sufficiently large, the results of Theorem \ref{t:hbex} hold and the direction of bifurcation is given by 
\beq
\mathrm{sign}[\theta_+] = -\,\mathrm{sign}\,\gamma. 
\eeq

\end{Theorem}
\begin{Proof}
Using Proposition \ref{p:efcns}, it is readily checked that Hypotheses \ref{h:l0}, \ref{h:HC}, and \ref{h:ex} are all satisfied. The existence of a Hopf bifurcation then follows by applying Theorem \ref{t:hbex}.

All that is left is to determine the sign of $\theta_+$.  In order to facilitate this determination, let $p$ and $\psi$ be as given in Proposition \ref{p:efcns} with normalization constants $A,B$ such that $A^3B = \re^{2\mu_\rlin \ell}$. Since the first order system vector $U_p$ decays exponentially fast outside the unstable interval, $(-\ell,\ell)$, the estimates \eqref{e:efcndec}, \eqref{e:efcndec2} in Proposition \ref{p:efcns} give 
\beq
\theta_+ =\int_{-\ell}^\ell \lp(3\p_u^3f(x,u_*(x)) p(x)^2 \overline{p(x)} \rp)_{xx} \overline{\psi(x)} +\mc{O}(\ell^{-4}).  
\eeq

The form of the solution $\tl W_{\rc,*}$ found in the proof of Proposition \ref{p:efcns} gives
\begin{align}\label{e:ABf}
\int_{-\ell}^\ell \lp(3\p_u^3f(x,u_*(x))  p(x)^2 \overline{p(x)} \rp)_{xx} \overline{\psi(x)} \,dx &=
18\gamma A^3 B \int_{-\ell}^\ell \re^{2\mu_\rlin x} \lp[(2\nu+ \overline{\nu})^2\sin^4\lp(\fr{\pi(x - \ell)}{2\ell}\rp) + \mc{O}(\ell^{-l})   \rp] dx\notag\\
&= -\fr{27\gamma(2\nu_\rlin+\overline{\nu_\rlin})^2}{8\mu_\rlin}+\mc{O}(\ell^{-1}).
 \end{align}
 Thus, for $\ell>0$ sufficiently large, 
 \begin{align}
\mathrm{sign}\lp[\,\rre\,\theta_+\rp]&=-\mathrm{sign}\lp[\,\rre\, \fr{27\gamma(2\nu_\rlin+\overline{\nu_\rlin})^2}{8\mu_\rlin}\rp]\notag\\
&= \mathrm{sign}\lp[\,\rre\, \gamma\,(2\nu_\rlin+\overline{\nu_\rlin})^2\rp] =\mathrm{sign}\lp[\gamma\,(  9\mu_\rlin^2 - \kappa_\rlin^2 )\rp]\notag\\
&= -\mathrm{sign}  \,\gamma
 \end{align}
 where the expressions given in Lemma \ref{l:calc} are used in the last two lines.

  \end{Proof}

\begin{Remark}
Since the argument of $\nu_\rlin$, and hence the sign of $9\mu_\rlin^2 - \kappa_\rlin^2$, is invariant with respect to changes in the value of $f'(\alpha)$ for $\alpha$ near zero, the sign of $9\mu_\rlin^2 - \kappa_\rlin^2$ will remain constant when our equation is linearized about a front $u_*\equiv \alpha$.
\end{Remark}

\begin{Remark}
We note that since $\lambda_\rlin$ is an accumulation point of the eigenvalues of $L$ as $l\rightarrow \infty$ successive Hopf bifurcations will rapidly occur as $c$ is decreased below $c_\rlin$.
\end{Remark}

\begin{Remark}
The findings of Theorem \ref{t:exam} are in agreement with numerical simulations, where supercritical behavior was found for $\gamma<0$ and subcritical behavior was found for $\gamma>0$. In the latter case, for $c$ slightly larger than $c_*$, we also observed hysteretic behavior between the front $u_*$ and a bifurcating periodic pattern.  This region of bistable, hysteretic behavior should, in principle, be able to be determined by finding higher order coefficients in the bifurcation equation.  We also note that the wavenumber of the periodic pattern was different than that predicted by the linearized equation, indicating that such solutions should be related to pushed fronts.
\end{Remark}

\section{Discussion}\label{s:dis}
Our methods should be applicable in many different settings. First of all, the existence result for viscous shocks in \cite{Sanstede08} can readily be obtained (and shortened significantly) with a nearly direct translation of our approach.  Also, problems with more general $u$-dependent source terms which are still exponentially localized in space could also be treated using our method. With such a source, the corresponding nonlinear solution operator $\mF$ would lose its conservation form. Since the codomain cannot be restricted as above, the method described in Remark \ref{r:par} must be employed to obtain a Fredholm index 0 operator. One such area where these source terms appear is in the equations governing the propagation of oscillatory detonation waves.  Here an ignition function, dependent on the characteristics of the gas, controls the reaction terms in the equation which feed the combustion; see \cite{Texier08}. These sources also arise in certain forms of the chemotaxis equation where the 
aggregation of bacteria depends nonlinearly on the density of bacteria (in addition to the gradient of the chemoattractant); see \cite{painter02}.  
 
Furthermore, our methods could be used to study problems with spatial dimension larger than one. In particular, for systems whose spatial domain is an infinite cylinder, Fredholm properties could be established using exponentially weighted spaces and a closed range lemma, while the index could be determined via a spectral flow.  This would then allow one to perform a Lyapunov-Schmidt reduction to obtain a bifurcation equation for transverse modes. Such an abstract functional analytic method will hopefully be simpler than the spatial dynamics formulations developed in \cite{peterhof97} and subsequent works, and the diffusive stability approach used by \cite{pogan14}.

\begin{Acknowledgment}
A. Scheel was partially supported by the National Science Foundation through grants DMS-0806614 and DMS-1311740. This material is based upon work supported by the National Science Foundation Graduate Research Fellowship under grant NSF-GFRP-00006595. Any opinion, findings, and conclusions or recommendations expressed in this material are those of the authors(s) and do not necessarily reflect the views of the National Science Foundation.
\end{Acknowledgment}

\bibliography{CH-ARMA}
\bibliographystyle{siam}
\end{document}